\documentclass[12pt]{article}
\usepackage{amssymb,amsmath,amstext,amsthm,amsfonts}
\usepackage[dvips]{graphicx}

\title{Statistical stability for robust classes of\\
maps with non-uniform expansion}
\author{Jos\'e F. Alves, Marcelo Viana
\thanks{The first author acknowledges the hospitality of IMPA,
where most of this work has been done, and financial support
from IMPA/CNPq and Funda\c c\~ao Calouste Gulbenkian.
The second author is partially supported by Pronex-Dynamical
Systems and Faperj. The authors are also greatful to the Paul Erd\H os Center in Budapest where the final ideas of this work were discussed.}
}

\date{\today}

\begin{document}

\newcommand{\vol}{\operatorname{vol }}
\newcommand{\graph}{\operatorname{graph }}

\def \RR {{\mathbb R}}
\def \ZZ {{\mathbb Z}}
\def \NN {{\mathbb N}}

\def \cc {{\cal C}}
\def \cad {{\cal D}}
\def \cf {{\cal F}}
\def \cl {{\cal L}}
\def \ch {{\cal H}}
\def \cn {{\cal N}}
\def \cp {{\cal P}}
\def \cq {{\cal Q}}
\def \car {{\cal R}}
\def \cs {{\cal S}}
\def \cau {{\cal U}}

\def \al {\alpha }
\def \be {\beta }
\def \si {\sigma }
\def \f {\varphi }
\def \F {\phi }
\def \t {\theta }
\def \th {\Theta }
\def \la {\lambda }
\def \La {\Lambda }
\def \o {\omega }
\def \O {\Omega }
\def \De {\Delta }
\def \de {\delta }
\def \ga {\gamma }
\def \e {\epsilon }

\newtheorem{maintheorem}{Theorem}
\renewcommand{\themaintheorem}{\Alph{maintheorem}}

\newtheorem{T}{Theorem}[section]
\newcommand{\ct}{\begin{T}}
\newcommand{\ft}{\end{T}}

\newtheorem{Corollary}[T]{Corollary}
\newcommand{\cco}{\begin{Corollary}}
\newcommand{\fco}{\end{Corollary}}

\newtheorem{Proposition}[T]{Proposition}
\newcommand{\cpr}{\begin{Proposition}}
\newcommand{\fpr}{\end{Proposition}}

\newtheorem{Lemma}[T]{Lemma}
\newcommand{\cle}{\begin{Lemma}}
\newcommand{\fle}{\end{Lemma}}

\newtheorem{Remark}[T]{Remark}
\newcommand{\cre}{\begin{Remark}}
\newcommand{\fre}{\end{Remark}}

\newtheorem{Definition}[T]{Definition}
\newcommand{\cd}{\begin{Definition}}
\newcommand{\fd}{\end{Definition}}

\newcommand{\dem}{\begin{proof}}
\newcommand{\cqd}{\end{proof}}

\newcommand{\var}{\mbox{var}\:}
\newcommand{\qand}{\quad\text{and}\quad}

\maketitle

\begin{abstract} We consider open sets of maps in a
manifold $M$ exhibiting non-uniform expanding behaviour in
some domain $S\subset M$. Assuming that there is a forward
invariant region containing $S$ where each map has a unique
SRB measure, we prove that under general uniformity
conditions, the SRB measure varies continuously in the
$L^1$-norm with the map.

As a main application we show that the open class of maps
introduced in [V] fits to this situation, thus proving that
the SRB measures constructed in [A] vary continuously with
the map. \end{abstract}

\section{Introduction}

In general terms, Dynamics has a twofold aim:
to describe, for the majority of dynamical systems, the typical behaviour
of trajectories, specially as time goes to infinity;
to understand how this behaviour changes when the system is modified, and
to what extent it is {\em stable\/} under small modifications.
In this work we are primarily concerned with the latter problem.

A first fundamental concept of stability, {\em structural stability\/},
was formulated by Andronov and Pontryagin~\cite{AP}.
It requires that the whole orbit structure remain unchanged under any
small perturbation of the dynamical system:
there exists a homeomorphism of the ambient manifold mapping trajectories
of the initial system onto trajectories of the perturbed one, preserving
the direction of time.
In the early sixties, Smale introduced the notion of {\em uniformly
hyperbolic\/} (or {\em Axiom A \/}) system, having as one of his main
goals to obtain a characterization of structural stability.
Such a characterization was conjectured by Palis and Smale in \cite{PS}:
a diffeomorphism (or a flow) is structurally stable if and only if it is
uniformly hyperbolic and satisfies the so-called strong transversality
condition.
Before that, structural stability had been proved for certain classes of
systems, including Anosov and Morse-Smale systems.
The ``if'' part of the conjecture was proved by Robbin, de Melo, Robinson
in the mid-seventies.
The converse remained a major open problem for yet another decade,
until it was settled by Ma\~n\'e for $C^1$ diffeomorphisms
(perturbations are small with respect to the $C^1$ norm).
The flow case was recently solved by Hayashi, also in the $C^1$ category
(the $C^k$ case, $k>1$, is still open both for diffeomorphisms and for flows).
See e.g. the book of Palis and Takens~\cite{PT} for precise definitions,
references and a detailed historical account.

Despite these remarkable successes, structural stability proved to be too
strong a requirement for many applications.
Several important models, including e.g. Lorenz flows and H\'enon maps,
are not stable in the structural sense, yet key aspects of their dynamical
behaviour clearly persist after small modifications of the system.
Weaker notions of stability, with a similar topological flavour,
were proposed throughout the sixties and the seventies,
but they all turned out to be too restrictive.

More recently, increasing emphasis has been put on expressing stability in
terms of persistence of {\em statistical properties\/} of the system.
A natural formulation, the one that concerns us most in this work,
corresponds to continuous variation of physical measures as a function
of the dynamical system. Let us explain this in precise terms.
We consider discrete-time systems, namely, smooth transformations
$\varphi:M\to M$ on a manifold $M$.
A Borel probability measure $\mu$ on $M$ is a {\em Sinai-Ruelle-Bowen
(SRB) measure\/} (or a {\em physical measure\/}), if there exists a
positive Lebesgue measure set of points $z\in M$ for which
\begin{equation}
\label{eq.SRB}
\lim_{n\to+\infty}\frac{1}{n}\sum_{j=0}^{n-1} f(\varphi^j(z))
=\int f\,d\mu
\end{equation}
for any continuous function $f:M\to\RR$.
In other words, time averages of all continuous functions are given by the
corresponding spatial averages computed with respect to $\mu$, at least
for a large set of initial states $z\in M$.

Let us suppose that $\varphi$ admits a forward invariant region
$U\subset M$, meaning that $\f(U)\subset U$, and there
exists a (unique) SRB measure $\mu=\mu_\varphi$ supported
in $U$ such that (\ref{eq.SRB}) holds for Lebesgue almost
every point $z\in U$. We say that $\varphi$ is {\em
statistically stable\/} (restricted to $U$) if similar
facts are true for any $C^k$ nearby map $\psi$, for some
$k\ge 1$, and the map $\psi\mapsto \mu_\psi$, associating
to each $\psi$ its SRB measure $\mu_\psi$, is continuous at
$\psi=\varphi$. For this definition, we consider in the
space of Borel measures the usual weak$^*$ topology: two
measures are close to each other if they assign close-by
integrals to each continuous function. Thus, this notion of
stability really means that time averages of continuous
functions are only slightly affected when the system is
perturbed.

Uniformly expanding smooth maps are well known to be
statistically stable, and so are Axiom A diffeomorphisms,
restricted to the basin of each attractor. On the other
hand, not much is known in this regard outside the
uniformly hyperbolic context. In the present work we
propose an approach to proving statistical stability for
certain robust (open) classes of non-uniformly expanding
maps. Precise conditions will be given in the next
subsection. For the time being, we just mention that
our
maps $\varphi$ exhibit asymptotic expansion,
 $$
\lim_{n\to+\infty}\frac{1}{n}\log\|D\f^n(z)v\| > 0
\quad\text{for every } v\in T_z M ,
 $$
 at Lebesgue almost
every point $z$ in some forward invariant region $U$ (but
they are not uniformly expanding). Moreover, they admit a
unique SRB measure which is an ergodic invariant measure
absolutely continuous with respect to Lebesgue measure in
$U$. These properties remain valid in a neighbourhood of
the initial map, and we prove that the SRB measure varies
continuously with the mapping in this neighbourhood. In
fact, our approach proves statistical stability in a strong
sense: the density (Radon-Nikodym derivative with respect
Lebesgue measure $m$) of the SRB measure,
${d\mu_\f}/{dm} $,
varies continuously with $\f$ as an
$L^1$-function.

To the best of our knowledge this is the first result of
statistical stability for maps with non-uniform expansion.
An application, and the example we had in mind when we
started this work, are the maps with multidimensional
non-uniform expansion introduced in \cite{V}, and whose SRB
measures were constructed in \cite{A}. In this context we
mention the important work of Dolgopyat~\cite{D}, where
statistical stability (and other ergodic properties) were
proved for some open classes of diffeomorphisms having
partially hyperbolic attractors whose central direction is
{\em mostly contracting\/} (negative Lyapunov exponents).
In that situation, cf. also Bonatti-Viana~\cite{BV}, SRB
measures are absolutely continuous with respect to Lebesgue
measure along the strong-unstable (uniformly expanding)
foliation of the attractor. Our systems in the present work
are closer in spirit to partially hyperbolic attractors
with {\em mostly expanding\/} central direction, in the
sense of Alves-Bonatti-Viana~\cite{ABV}. Statistical
stability for the latter systems has not yet been proved.


\subsection{Statement of results}\label{s.statement}

Let $\f \colon M \rightarrow M$ be a map from some $d-$dimensional
manifold into itself, $S$ be some region in $M$, and $\F:S\rightarrow S$
be a return map for $\f$ in $S$.
That is, there exists a countable partition $\car=\{R_i\}_i$
into subsets of $S$,
and there exists a function $h\colon\car\rightarrow\ZZ^+$ such that
$$
\F\mid R=\f^{h(R)}\mid R
\quad\mbox{ for each }\quad R\in \car.
$$
For simplicity, we will assume that $S$ is diffeomorphic to some
bounded region $\tilde{S}$ of $\RR^n$ (but similar arguments hold
in general, using local charts).
Then we can pretend that $S\subset\RR^n$, through identifying it
with $\tilde{S}$, and we do so.

We say that $\F$ is a {\em $C^2$ piecewise expanding map\/}
if the following conditions hold:
 \begin{itemize}
\item[(1)] The boundary of each $R_i$ is piecewise $C^2$ (a
countable union of $C^2$ hyper-surfaces) and has finite
$(d-1)$-dimensional volume.
\item[(2)]  Each $\F_i\equiv\F\mid R_i$ is a  $C^2$
bijection from the interior of $R_i$ onto its image,
admitting a $C^2$ extension to the closure of $R_i$.
\item[(3)]  There is $0<\si<1$ such that $\| D\F_i^{-1}\|
<\si$ for every $i\geq 1$.
 \end{itemize} We say that $\F$
satisfies a {\em bounded distortion\/} property if:
\begin{itemize} \item[(4)] There is some $K>0$ such that
for every $i\geq 1$
 $$
 \frac{\left\|
D\left(J\circ\F^{-1}_i\right)\right\|}
{\left|J\circ\F^{-1}_i\right|}<K,
 $$
 where $J$ is the Jacobian of $\F$.
  \end{itemize}
 Moreover, we assume that
the images under $\F$ of all the elements of the partition
$\car$ satisfy the following {\em bounded geometry\/}
condition:
 \begin{itemize}
 \item[(5)] There are constants
$1\geq \be>\si/(1-\si)$ and $\rho>0$ such that the boundary
of each $R_i$ has a tubular neighborhood of size $\rho$
inside $R_i$, and the $C^2$ components of the boundary of
each $R_i$ meet at angles greater than $\arcsin(\be)>0$.
 \end{itemize}


It was shown in \cite[Section 5]{A} that conditions (1)--(5)
imply that the map $\F$ has some invariant
probability measure $\mu$ absolutely continuous with respect to
Lebesgue measure on $S$ (henceforth denoted $m$ and assumed
to be normalized).
Then
$$
\mu^{\ast} =\sum_{j=0}^{\infty}\f_{\ast}^j\left(\mu\mid
\{h>j\}\right)
 $$
is an absolutely continuous invariant measure for $\f$.
Moreover, the density $d\mu/dm$ of $\mu$ is in $L^p(S)$ for
$p=d/d-1$. As a consequence, the measure $\mu^{\ast}$ is
finite, as long as we have:
 \begin{itemize} \item[(6)]
  The function $h$ is in $L^q(S)$ for $q=d$\quad
          (this is taken so that $1/p+1/q=1$).
 \end{itemize}
It was also observed in \cite[Sections 5 and 6]{A} that the
absolutely continuous invariant measure $\mu^{\ast}$ may be
taken ergodic (which implies that it is an SRB measure for
$\f$) and, moreover, $\f$ has finitely many such ergodic
measures.


\medskip

Now we state our first main result. Let $k\ge 1$ be fixed,
and $\cau$ be an open set of $C^k$ transformations on $M$
admitting a forward invariant compact region $U$. We endow $\cau$
with the $C^k$ topology. Assume that we may associate to
each $\f\in\cau$ a map $\F_\f\colon S\rightarrow S$, a
partition $\car_\f$ of $S\subset U$, and a function
$h_\f\colon\car_\f\rightarrow \ZZ^+$, satisfying properties
(1) to (6) above. We consider elements $\f_0$ of $\cau$
satisfying the following uniformity conditions:
 \begin{itemize}
 \item[(U1)] Given
any integer $N\geq 1$ and any $\e>0$, there is
$\de=\de(\e,N)>0$ such that for $j=1,\dots,N$
 $$
 \|\f-\f_0\|_{C^k}<\de\Rightarrow
m\big(\{h_\f=j\}\Delta\{h_{\f_0}=j\}\big)<\e,
 $$
where $\Delta$ represents symmetric difference of two sets.
\item[(U2)] Given $\e>0$, there are $N\geq 1$ and $\de=\de(\e,N)>0$
 for which
 $$
 \|\f-\f_0\|_{C^k}<\de\Rightarrow
 \big\|\sum_{j=N}^\infty\chi_{\{h_\f>j\}}\big\|_q<\e.
 $$

\item[(U3)] Constants $\sigma$, $K$, $\beta$, $\rho$ as
above may be chosen uniformly in a $C^k$ neighborhood of
$\f_0$\,.
 \end{itemize}

\begin{maintheorem}
\label{t.A}
Let $\cau$ be as above, and suppose that every $\f\in\cau$ admits a
unique SRB measure $\mu_\f$ in $S$.
Then
\begin{enumerate}
\item $\mu_\f$ is absolutely continuous with respect to the Lebesgue
      measure $m$;
\item if $\f_0\in\cau$ satisfies (U1), (U2), (U3) then $\f_0$ is
      statistically stable in a strong sense: the map
 $$
 \cau\ni\f \mapsto \frac{d\mu_\f}{dm}
 $$
 is continuous, with respect to the $L^1$-norm, at $\f=\f_0$.
\end{enumerate} \end{maintheorem}



We observe that under assumption (U1), condition (U2) can
be reformulated in equivalent terms as:
 \begin{itemize}
 \item[(U2')] Given $\e>0$, there is  $\de>0$ for which
 $$
\|\f-\f_0\|_{C^k}<\de\Rightarrow \|h_\f-h_{\f_0}\|_q<\e.
 $$
 \end{itemize}
 A simple proof of this equivalence will be given in  Section \ref{s.stability} (just before Proposition \ref{p.cont*}).

\medskip

Our next results state that the assumptions of Theorem~\ref{t.A} do
correspond to robust classes of smooth maps in some manifolds.

\begin{maintheorem} \label{t.B} There exists a non-empty
open set $\cn$ in the space of $C^3$ transformations from
$S^1\times I$ into itself such that conditions (1)--(6) and
(U1)--(U3) are satisfied by every element of $\cn$.
\end{maintheorem}

The open set $\cn$ we exhibit for the proof of this result is the
one constructed in \cite{V}.
As pointed out in that paper, the choice of the cylinder $S^1\times I$,
$I=[0,1]$, as ambient space is rather arbitrary, the construction
extends easily to more general manifolds.
In what follows we briefly describe the set $\cn$, referring the
reader to \cite{V} and Section~\ref{s.hyperbolicreturns} for more
details.

Let $d$ be some large integer: $d\ge 16$ suffices, but this
is far from being optimal. Let $a_0\in(1,2)$ be such that
the critical point $x=0$ is pre-periodic under iteration by
the quadratic map $q(x)=a_0-x^2$ (again, this is far too
strong a requirement on the parameter $a_0$). Let
$b:S^1\rightarrow \RR$ be a Morse function, for instance,
$b(t)=\sin(2\pi t)$. Note that $S^1=\RR/\ZZ$. For each
$\alpha>0$, consider the map $\f_{\al}:S^1\times
\RR\rightarrow S^1\times \RR$ given by
$\f_\al(\t,x)=(\hat{g}(\t), \hat f(\t,x))$, where $\hat{g}$
is the uniformly expanding map of the circle defined by
$\hat{g}(\theta)=d\theta$ (mod $\ZZ$), and $\hat
f(\t,x)=a(\t)-x^2$ with $a(\t)=a_0+\al b(\t)$. We shall
take $\cn$ to be a small $C^3$ neighborhood of $\f_\alpha$,
for some (fixed) sufficiently small $\alpha$.

It is easy to check that for $\al$ small enough there is an
interval $I\subset (-2,2)$ for which $\f_\al(S^1\times I)$
is contained in the interior of $S^1\times I$. Thus, any
map $\f$ close to $\f_\al$ in the $C^0$ topology has
$U=S^1\times I$ as a forward invariant region, and so $\f$
has an attractor inside this invariant region,
which is precisely the set
$$
\Lambda=\bigcap_{n\geq 0}\f^n(U).
$$


As we mentioned before, properties (1)--(6) imply that the maps
$\f\in\cn$ admit (finitely many) SRB measures, which are ergodic
absolutely continuous invariant measures.
In order to be able to apply Theorem~\ref{t.A} to this open set
$\cn$, we also have to show that the SRB measure is unique for
each $\f\in\cn$.
This will follow from a stronger fact, stated in Theorem~\ref{t.C}
below.

Let us say that $\f$ is {\em topologically mixing} if for every
open set $A\subset S^1\times I$ there is some $n=n(A)\in\ZZ^+$
for which $\f^n(A)=\Lambda$, and say that $\f$ is {\em ergodic
with respect to Lebesgue measure} if for every Borel subset
$B\subset S^1\times I$ such that $\f^{-1}(B)=B$, either $B$ or
$(S^1\times I)\setminus B$ have Lebesgue measure equal to zero.

\begin{maintheorem}
\label{t.C}
Let $\cn$ be as described above.
Then the transformations $\f\in\cn$ are topologically mixing
and ergodic with respect to Lebesgue measure.
\end{maintheorem}


\section{Absolute continuity}\label{s.absolute}

A main ingredient in our arguments, as well as in \cite{A},
is the notion of variation for functions in higher
dimensions. For $f\in L^1(\RR^d)$ with compact support we
define the {\it variation} of $f$ as $$
\var(f)=\sup\left\{\int_{\RR^d}f\mbox{div}(g)dm\colon g\in
C_0^1(\RR^d,\RR^d), \|g\|_0\leq 1\right\}, $$ where
$C_0^1(\RR^d,\RR^d)$ is the set of $C^1$ maps from $\RR^d$
to $\RR^d$ with compact support and $\|\:\:\|_0$ is the
supremum norm in $C_0^1(\RR^d,\RR^d)$. We observe that in
the case of $f$ be a $C^1$ map, then $\var(f)$ coincides
with $\int \|Df\| dm$ (see e.g. \cite[Example 1.2]{G}). We
consider the space of {\em bounded variation} functions
 $$
BV(\RR^d)=\left\{ f\in
L^1(\RR^d)\colon\var(f)<+\infty\right\}.
 $$
 We will make
use of the following results concerning bounded variation
functions:

\cpr \label{p.approx}
 Given $f\in BV(\RR^d)$, there is a sequence $(f_n)_n$ of $C^\infty$
 maps such that
 $$
 \lim_{n\rightarrow\infty}\int|f-f_n|dm=0
 \qand
 \lim_{n\rightarrow\infty}\int\|Df_n\|dm=var(f).
 $$
 \fpr
 \dem See \cite[Theorem 1.17]{G}. \cqd

\cpr \label{p.compact}
If $(f_k)_k$ is a sequence of functions in $BV(\RR^d)$ such that there
is a constant $K_0>0$ for which
$$
\var(f_k)\leq K_0\quad\mbox{and}\quad
\int|f_k|dm\leq K_0\quad\mbox{for every }k,
$$
then $(f_k)_k$ has a subsequence converging in the $L^1$-norm to an
$f_0$ with $\var(f_0)\leq K_0$.
\fpr
\dem See \cite[Theorem 1.19]{G}. \cqd

\cpr \label{p.dimension}
 Let $f\in BV(\RR^d)$ and take $p=d/(d-1)$. Then
 $$
 \|f\|_p\leq K_1\var(f),
 $$
 where $K_1>0$ is a constant depending only on $d$.
 \fpr \dem See \cite[Theorem 1.28]{G}. \cqd

Now we introduce the linear transfer operator associated to
$\F$, $$ \cl_\F:L^1(S)\longrightarrow L^1(S) $$ defined as
$$ \cl_\F f=\sum_{i=1}^{\infty}\frac{f\circ\F_i^{-1}}
{|J\circ\F_i^{-1}|}\chi_{\F(R_i)}. $$ It is well-known that
each fixed point of $\cl_\F$ is the density of an
absolutely continuous $\F$-invariant finite measure. The
next lemma gives a Lasota-Yorke type inequality for maps in
$BV(\RR^d)$, which plays a crucial role in the proof of the
existence of fixed points for $\cl_\F$.

\cle \label{l.ly} There are constants $0<\la<1$ and $K_2>0$
such that for every $f\in BV(\RR^d)$,
 $$
 \var(\cl_\F f)\leq\la\var(f)+K_2\int |f|dm.
 $$ \fle \dem See \cite[Lemma 5.4]{A}. \cqd

\cre\label{r.unif} It follows from the proof of \cite[Lemma
5.4]{A} that the constant $\la$ is equal to
$\si(1+1/\beta)$. Hence, by assumption (U3), $\la$ may be
taken uniformly smaller than one in a whole neighborhood of
$\f$. It also follows from the proof of \cite[Lemma 5.4]{A}
that the constant $K_2$ coincides with
$K+1/(\beta\rho)+K\beta$, which may also be taken uniform
in a neighborhood of the map $\f$. Hence, the constants
$0<\la<1$ and $K_2>0$ may be taken in such a way that the
Lasota-Yorke type inequality in the previous lemma holds
for every map in a neighborhood of $\f\in\cau$. \fre

Consider for each $k\geq 1$ the function
 $$
f_k=\frac{1}{k}\sum_{j=0}^{k-1}\cl_\F^j 1.
 $$
  We have
   \[
\int |f_k|dm=1\quad\mbox{for }k\geq 1,
 \]
 and it follows
from Lemma \ref{l.ly} that
 \[ \var(f_k)\leq
K_3\quad\mbox{for }k\geq 1,
 \] where
$K_3=\var(\chi_S)+K_2\sum^\infty_{k=0}\la^k+1.$ It follows
from Proposition \ref{p.compact} that $(f_k)_k$ has a
subsequence converging in the $L^1$-norm to some $\rho$
with $\var(\rho)\leq K_3$. Hence, $\mu_\F=\rho m$ is an
absolutely continuous $\F$-invariant probability measure.
From this it is deduced in \cite[Section 6]{A} that
 $$
\mu^{\ast}_\f
=\sum_{j=0}^{\infty}\f_{\ast}^j\left(\mu_\F\mid
\{h_\f>j\}\right)
 $$
is an absolutely continuous $\f$-invariant
finite measure.

 \cle\label{l.sup}
 Given any $\F$-invariant set $A\subset S$ with positive
 Lebesgue measure, there is an absolutely continuous
 $\F$-invariant probability measure $\mu_A=f_Am$ for which
 $\mu_A(A)=1$. Moreover, $f_A$  may be taken
 with $\var (f_A)\leq 4K_2$.
\fle
 \dem
 We start by proving that given any
$f\in L^1(\RR^d)$, the sequence
$1/n\sum_{j=0}^{n-1}\cl_\F^j f$ has accumulation points (in
the $L^1$-norm) in $BV(\RR^d)$. Let $f\in L^1(\RR^d)$ and
take a sequence $(f_n)_n$ in $BV(\RR^d)$ converging to $f$
in the $L^1$-norm. It is no restriction to assume that
$\|f_n\|_1\leq 2\|f\|_1$ for every $n\geq 1$ and we do it.
For each $n\geq 1$ we have
 $$
 \var(\cl^j_\F f_n)\leq \lambda^j\var(f_n)+K_2\|f_n\|_1\leq 3K_2\|f\|_1
 $$
for  large $j$. So, taking $k$ large enough we have
 $$
\var\left(\frac{1}{k}\sum_{j=0}^{k-1}\cl^j_\F
f_n\right)\leq 4K_1\|f\|_1.
 $$
 Using the well known fact that the transfer operator does
 not expand $L^1$-norms, we also have
 $$
 \big\|\frac{1}{k}\sum_{j=0}^{k-1}\cl^j_\F f_n\big\|_1\leq
 \frac{1}{k}\sum_{j=0}^{k-1}\|\cl^j_\F f_n\|_1\leq 2\|f\|_1
 $$
 for every  $j\geq 1$.
 It follows from Proposition \ref{p.compact} that there exists some
$\hat f_n\in BV(\RR^d)$ and a sequence $(k_i)_i$ for which
 $$
 \lim_{i\rightarrow\infty}
 \big\|\frac{1}{k_i}\sum_{j=0}^{k_i-1}\cl^j_\F f_n-\hat f_n\big\|_1=0
 $$
 and, moreover,  $\var(\hat f_n)\leq 4K_2\|f\|_1$. Now we
apply the same argument to the sequence $(\hat f_n)_n$ in
order to obtain a subsequence $(n_l)_l$ such that $(\hat
f_{n_l})_l$ converges in the $L^1$-norm to some $\hat f$
with $\var(\hat f)\leq 4K_2\|f\|_1$. Since
 $$
 \big\|\frac{1}{k}\sum_{j=0}^{k-1}\cl^j_\F f_{n_l}-\hat
 f\big\|_1\leq
 \big\|\frac{1}{k}\sum_{j=0}^{k-1}\cl^j_\F f_{n_l}-\hat
f_{n_l}\big\|_1+\|\hat f_{n_l}-\hat f\|_1,
 $$
 there is some sequence $(k_l)_l$ for which
 $$
 \lim_{l\rightarrow\infty}
 \big\|\frac{1}{k_l}\sum_{j=0}^{k_l-1}\cl^j_\F f_{n_l}-\hat
f\big\|_1=0.
 $$
On the other hand,
 $$
\big\|\frac{1}{k_l}\sum_{j=0}^{k_l-1}\big(\cl^j_\F
f_{n_l}-\cl^j_\F
f\big)\big\|_1\leq\frac{1}{k_l}\sum_{j=0}^{k_l-1}
\left\|f_{n_l}-f\right\|_1=\left\|f_{n_l}-f\right\|_1
 $$ and
this last term goes to 0 as $l\rightarrow \infty$. This
finally implies that
 $$\lim_{l\rightarrow\infty}
 \big\|\frac{1}{k_l}\sum_{j=0}^{k_l-1}\cl^j_\F f-\hat
 f\big\|_1=0,
 $$
thus proving that the sequence $1/n\sum_{j=0}^{n-1}\cl^j_\F
f$ has accumulation points in $BV(\RR^d)$.

Now let $A\subset S$ be a $\phi$-invariant set with
positive Lebesgue measure. Considering $f=\chi_A\in
L^1(\RR^d)$ in the previous argument, we obtain some
sequence $(k_l)_l$ and $f_A\in BV(\RR^d)$ for which
 $$
  \lim_{l\rightarrow\infty}
 \big\|\frac{1}{k_l}\sum_{j=0}^{k_l-1}\cl^j_\F\chi_A-f_A\big\|_1=0
 $$
 Moreover,
 $$
\var(f_A)\leq 4K_2\|\chi_A\|_1\leq 4K_2\qand\|f_A\|_1>0
 $$
 (here we use that $m\mid S$ is normalized). Considering
 $f_A$ already multiplied by a constant factor in order to have
  $L^1$-norm equal to 1,
 we take $\mu_A=f_Am$. Since
 $$
\mu_A(S\setminus A)=
\lim_{l\rightarrow\infty}\frac{1}{k_l}\sum_{j=0}^{k_l-1}
\int_{S\setminus A}\cl^j_\F\chi_A\;dm
=\lim_{l\rightarrow\infty}\frac{1}{k_l}\sum_{j=0}^{k_l-1}
\int_S\chi_{(S\setminus A)}\circ\phi^j\cdot \chi_A\;dm = 0
 $$
 we have that $\mu_A$ gives full weight to $A$, thus
concluding the proof of the result. \cqd

\cco\label{c.inf}
 There is a constant $\widehat K(d)>0$ such that if
 $A\subset S$ is a $\F$-invariant set
 with positive Lebesgue measure, then
 $m(A)\geq \widehat K(d)$.
\fco
 \dem
Let $A\subset S$ be a $\F$-invariant set
 with positive Lebesgue measure and $\mu_A=f_Am$ a measure as in
 Lemma \ref{l.sup}.  Since $f_A\in BV(\RR^d)\subset
 L^p(\RR^d)$ (recall Proposition \ref{p.dimension})
 and $\mu_A$ gives full weight to $A$, it follows from Minkowski's
 inequality that
 \[
  1=\int f_A dm\leq \|f_A\|_p\cdot\|\chi_A\|_q\leq K_14K_2\;m(A)^{1/d}.
 \]
 We take
$\widehat{K}(d)=(K_14K_2)^{-d}$.
\cqd

The  proposition below gives the first item of Theorem
\ref{t.A}.

\cpr\label{p.abscont}
 If $\f\in\cau$ has a unique SRB
measure $\mu_\f$ in $U$, then $\mu_\f^*=\mu_\f$.
 \fpr
 \dem
We will prove that $\mu_\f^*$ is ergodic, and so an SRB
measure for $\f$ in $U$. Let $A\subset U$ be any
$\f$-invariant set with $m(A\cap S)>0$. We have
 $$
\phi^{-1}(A\cap S) = \big\{x\in S \colon\phi(x)\in A\big\}
=\bigcup_{j\geq 1} \left(\f^{-j}(A)\cap \{h_\f=j\} \right)
 =A\cap S,
 $$
and so the set  $A\cap S$ is also $\phi$-invariant. Since
we are taking $A$ with $m(A\cap S)>0$, it follows from
Corollary \ref{c.inf} that $m(A\cap S)\geq \widehat{K}(d)$.
As a consequence, $S$ can be covered by a finite number of
$\f$-invariant sets $A_1,\dots,A_r$ intersecting $S$ in a
positive Lebesgue measure set, and which are minimal: for
$1\leq i\leq r$ there is no $\f$-invariant set $B_i\subset
A_i$ with $m(B_i\cap S)>0$.

It follows from Lemma \ref{l.sup} that for each
$i=1,\dots,r$ there is an absolutely continuous
$\F$-invariant measure $\mu_i$ giving full weight to
$A_i\cap S$. Take
 $$
\mu_i^{\ast} =\sum_{j=0}^{\infty}\f_{\ast}^j\left(\mu_i\mid
\{h_\f>j\}\right)
 $$
and denote $A_i^c=U\setminus A_i$. Since $A_i$ is
$\f$-invariant, we have that $A_i^c$ is also
$\f$-invariant, and so
 $$
\mu_i^{\ast}(A_i^c) =
\sum_{j=0}^{\infty}\mu_i\left(\f^{-j}\left(A_i^c\right)\cap
\{h_\f>j\}\right)= \sum_{j=0}^{\infty}\mu_i\left(A_i^c\cap
\{h_\f>j\}\right)=0
 $$
Hence, assuming  $\mu_i^\ast$  normalized we have that each
$\mu_i^\ast$ is  an absolutely continuous $\f$-invariant
probability measure giving full weight to $A_i$. The
minimality of each $A_i$ implies that $\mu_i^*$ is ergodic
for $1\leq i\leq r$, and so it coincides with the SRB
measure $\mu_\f$. This in particular implies that $r=1$.
The fact that $A_1$ is a minimal $\f$-invariant set that
contains $S$ implies that $\mu_\f^*$ is an ergodic
absolutely continuous $\f$-invariant probability measure,
thus coinciding with $\mu_\f$.
\cqd


\cre\label{r.decompor}
 The proof of Proposition \ref{p.abscont} gives also that
 under the hypothesis of $\f$ having a unique SRB
 measure in $U$, the region $S$  intersects a unique
 $\f$-invariant minimal set and, consequently, is contained in it.
 However, if we do not assume uniqueness of the SRB measure
 in $U$, we may write
 $\mu_\f^*=\sum_i\mu_\f^*(A_i)\mu^*_i$, where
 the sum is over the values of $i$ for which $\mu_\f^*(A_i)>0$
($A_1,\cdots A_r$ are the minimal sets given by
 the proof of Proposition \ref{p.abscont}) and
 each $\mu_i^*$ is the  normalized
 restriction of $\mu_\f^*$ to $A_i$, thus an SRB measure.
 More generally, if $\mu$ is an absolutely continuous
 $\f$-invariant probability measure giving full weight to
 $A_1\cup\cdots\cup A_r$, then $\mu=\sum_i\mu(A_i)\mu_i$
 where each $\mu_i$ is an SRB measure defined in the same
 way as $\mu^*_i$ above.
 \fre

\section{Statistical stability}\label{s.stability}

Now we prove that under the assumptions of
Theorem~\ref{t.A} the density of the measure $\mu^*_\f$
varies continuously in the $L^1$-norm with the map $\f$.
Let $\f_0$ be some map in $\cau$ and $(\f_n)_n$ a sequence
of maps in $\cau$ converging to $\f_0$ in the $C^k$
topology. As described above, we make for each on the maps
$\f_n$ and $\f_0$ the following choices:
 \[
\begin{array}{ccccccccccl}
 \f_n&\longmapsto&h_n &\longmapsto
&\F_n&\longmapsto&\rho_n&\longmapsto&\mu_n&\longmapsto&\mu^*_n
\\ \downarrow &&&&&&&&&& \downarrow \:\:?\\
\f_0&\longmapsto&h_0&\longmapsto&\F_0&\longmapsto&\rho_0&
\longmapsto&\mu_0&\longmapsto&\mu^*_0
 \end{array}
 \]
 It follows from the way we obtain each $\rho_n$ that
 $$
\var(\rho_n)\leq K_3\quad\mbox{and}\quad \int\rho_ndm\leq 1
 $$
for every $n\geq 1$ (recall also Remark \ref{r.unif}).
Thus, we may apply Proposition \ref{p.compact} to the
sequence of densities $(\rho_n)_n$ and deduce that it has
some subsequence $(\rho_{n_i})_{n_i}$ converging  in the
$L^1$-norm to some $\rho_\infty$ with
$\var(\rho_\infty)\leq K_3$. We consider
$\mu_\infty=\rho_\infty m$ and define
 $$ \mu^{\ast}_\infty =\sum_{j=0}^{\infty}\f_*^j\left(\mu_\infty\mid
 \{h_0>j\}\right).
 $$
The goal of the results below is to show that the densities
of $\mu^*_n$ with respect to the Lebesgue measure converge
in the $L^1$-norm to the density of $\mu_\infty^*$ and,
moreover, the measure $\mu_\infty^*$ coincides with $\mu^*_0$.
We start with some auxiliary lemmas.

\cle\label{l.variation}
There is  $K_4>0$ (depending only
on the dimension $d$) such that for every  $f\in BV(\RR^d)$
and $\psi\colon D\rightarrow\RR^d$ a $C^1$ diffeomorphism
from a compact $D\subset\RR^d$ onto its image
 $$
 \int_D|f\circ \psi-f|dm\leq K_4\|\psi-id\|_0^d\var(f).
 $$
 \fle
 \dem
We start by proving the result in the case of $f$ being
a continuous piecewise affine map, i.e. letting $\De$ be
the support of $f$, there is a finite number of domains
$\De_1,\dots,\De_N$ for which $\De=\cup_{i=1}^N \De_i$ and
$\nabla f$ (the gradient of $f$) is constant on each
$\Delta_i$.  We define
 $
 D_1=\{ (x,z)\in \RR^{d+1}\colon x\in D\qand z\in [f(x),f(\psi(x))]\}
 $
and $D_2$ the horizontal $\|\psi-id\|_{0}$-neighborhood of
the graph of $f$. That is,
 $D_2$ is equal to the set of points $(x,z)\in \RR^{d+1}$ for which
 there is $t\in\RR^d$ with $\|t\|\leq\|\psi-id\|_{0}$ and
$y\in \RR^d$ such that $x=y+t$ and $z=f(y)$. We observe
that $D_1\subset D_2$. Indeed, given $(x,z)\in D_1$, and
since $z\in [f(x),f(\psi(x))]$, by the continuity of $f$
there is $y\in [x,\psi(x)]$ such that $z=f(y)$. Taking
$t=y-x$ we have $\|t\|\leq \|\psi(x)-x\|$. Hence
 $$
\int_D|f\circ
\psi-f|dm=\int_D\int_{[f(x),f(\psi(x))]}1\:dzdm(x)=
\vol(D_1)\leq \vol(D_2).
 $$
For each $i=1,\cdots,N$ we define $G_i=\graph(f\mid \De_i)$
and $H_i$ the horizontal $\|\psi-id\|_0$-neighborhood of
$G_i$. We have
 $$
\vol(D_2)\leq \sum_{i=1}^N \vol(H_i).
 $$
Letting $\nabla_i f$ denote the gradient (constant) vector
of $f\mid \De_i$, we have that $(-\nabla_i f,1)$ is
orthogonal to $G_i$. Taking
$\partial_z=(0,\dots,0,1)\in\RR^{d+1}$ we have
 $$
\vol(H_i)\leq K_4\|\psi-id\|_0^d\sin\angle \big( (-\nabla_i
f,1),\partial_z\big)\vol(G_i),
 $$
where $K_4>0$ is a constant depending only on the volume of
the unit ball in $\RR^d$. We have
 $$
 \sin\angle \big((-\nabla_i f,1),\partial_z\big)=
 \sqrt{1-\cos^2\angle \big(
(-\nabla_i f,1),\partial_z\big)}= \frac{\|\nabla_i f\|}
{\sqrt{1+\|\nabla_i f\|^2}}
 $$
 and
  $$
 \vol(G_j)=\sqrt{1+\|\nabla_if\|^2}\vol(\De_i).
 $$
Altogether this yields
 $$
\int_D|f\circ \psi-f|dm\leq  K_4\|\psi-id\|_0^d
\sum_{i=1}^N \|\nabla_if\|\vol(\De_i).
 $$
 Taking into account that in this case
 $$
 \sum_{i=1}^N
\|\nabla_if\|\vol(\De_i)=\int\|\nabla f\|dm=\var(f),
 $$
 we obtain the result for any continuous piecewise affine map.

 The next step is to deduce the result for any
 $C^1$ map $f$. In this case, we may take a
 sequence $(f_n)_n$ of continuous piecewise affine maps
 such that
 $$
 \|f-f_n\|_0\rightarrow 0\qand \|Df-Df_n\|_0\rightarrow 0
 \quad\mbox{as}\quad n\rightarrow\infty
 $$
 (here we are take derivatives only in the
 interior points of the elements of the partition
 associated to each piecewise affine map). This implies
 that
 $$
 \int_D|f\circ\psi-f|dm=\lim_{n\rightarrow\infty}
 \int_D|f_n\circ\psi-f_n|dm
 $$
 and
 $$
 \var(f)=\int\|Df\|dm=\lim_{n\rightarrow\infty}\int\|Df_n\|dm
 =\lim_{n\rightarrow\infty}\var(f_n).
 $$
 Taking into account the case we have seen before, this
 implies the result also for the case of $f$ being a $C^1$ map.

 For the general case, we know by Proposition \ref{p.approx}
 that given $f\in BV(\RR^d)$ there is a sequence $(f_n)_n$
 of $C^1$ maps for which
 \begin{equation}\label{e.approx}
 \lim_{n\rightarrow\infty}\int|f-f_n|dm=0
 \qand
 \lim_{n\rightarrow\infty}\var(f_n)=\var(f).
 \end{equation}
 We have
 $$
 \int_D
 |f\circ\psi-f|dm\leq\int_D|f\circ\psi-f_n\circ\psi|dm
 +\int_D|f_n\circ\psi-f_n|dm
 +\int_D|f_n-f|dm.
 $$
 Since
 $$
 \int_D|f_n\circ
 \psi-f\circ\psi|dm=\int_{\psi(D)}|f_n-f|\cdot\Psi dm
 \leq \|\Psi\|_0\int|f_n-f|dm
 $$
 where $\Psi=1/|\det D\psi|\circ \psi^{-1}$, the result
 for general $f\in BV(\RR^d)$ follows from
 (\ref{e.approx}) and the previous case.
\cqd

At this point we also introduce the transfer operator
$\cl_\f$ associated to $\f\in\cau$, defined for each $f\in
L^1(\RR^d)$  as
 \begin{equation}\label{e.transfer}
 \cl_\f f(y)=\sum_{x\in \f^{-1}(y)}
\frac{f(x)}{|\det D\f(x)|}.
 \end{equation}
For our purposes the value of $\cl_\f f(y)$ for $y$ a
critical value of $\f$ is rather unimportant.
$\cl_\f$ is defined in such a way that
 $$
 \int (\cl_\f f)g\: dm=\int f(g\circ \f) \: dm
 $$
 for every $f,g\in L^1(\RR^d)$ wherever these integrals
 make sense. Let us say that $\cl_\f$ is being introduced
 just for the sake of notational simplicity, and so we stay
 away from rigorous formalities.

\cle\label{l.converge}
 Given $\e>0$ there is $\de>0$ such that
if  $\|\f-\f_0\|_{C^1}<\de$, then for every $f\in
BV(\RR^d)$ with support contained in $S$
 $$
 \int|\cl_\f f-\cl_{\f_0}f|dm<\e\var(f).
 $$
 \fle
\dem
 Take some small $\e_1>0$ and define $\cc(\e_1)$ the
$\e_1$-neighborhood of the critical set of $\f_0$ in $S$.
We divide $S\setminus\cc(\e_1)$ into a finite number of
domains of injectivity of $\f_0$ whose collection we call
$\cad(\f_0)$. We observe that if $\f$ is close enough to
$\f_0$, then $\cc(\e_1)$ also contains the critical set of
$\f$, and so we may define an analogous $\cad(\f)$ for $\f$
in such a way that for each $D_0\in\cad(\f_0)$ there is one
(and only one) $D\in\cad(\f)$ for which the Lebesgue
measure of $D\Delta D_0$ is small. For each
$D_0\in\cad(\f_0)$ let $D$ be the element in $\cad(\f_0)$
naturally associated to $D_0$, and define
 $$
\widehat{D}_0=\f_0^{-1}\big(\f_0(D_0)\cap \f(D)\big)=
D_0\cap\f_0^{-1}\circ\f(D).
 $$
 We have
  \begin{eqnarray}
\int |\cl_\f f-\cl_{\f_0}f |dm &\leq &
\int_{\f_0(\cc(\e_1))}|\cl_\f f-\cl_{\f_0}f |dm
\label{e1}\\ && +\sum_{D_0\in\cad(\f_0)}\int_{\f_0(D_0)\cap
\f(D)}|\cl_\f f-\cl_{\f_0}f|dm\label{e2}\\ &&
+\sum_{D_0\in\cad(\f_0)}\int_{\f_0(D_0)\setminus\f(D)}|\cl_{\f_0}f|dm
\label{e3}\\ &&
+\sum_{D_0\in\cad(\f_0)}\int_{\f(D)\setminus
\f_0(D_0)}|\cl_\f f|dm \label{e4} \end{eqnarray} Let us
estimate the expressions on the right hand side of the
inequality above. For the first one we have
 \begin{eqnarray*}
 \int_{\f_0(\cc(\e_1))}|\cl_\f
f-\cl_{\f_0}f |dm &\leq & \int_{\f_0(\cc(\e_1))}|\cl_\f
f|dm +\int_{\f_0(\cc(\e_1))}|\cl_{\f_0}f |dm \\ &\leq &
\int\chi_{\f_0(\cc(\e_1))}\cl_{\f}|f|dm
+\int\chi_{\f_0(\cc(\e_1))}\cl_{\f_0}|f|dm \\ &\leq &
\int\chi_{\f_0(\cc(\e_1))}\circ\f|f| dm
+\int\chi_{\f_0(\cc(\e_1))}\circ\f_0|f|dm
 \end{eqnarray*}
Since $f\in L^p(\RR^d)$
(recall Proposition \ref{p.dimension}), it
follows from Minkowski's inequality that
 $$
\int\chi_{\f_0(\cc(\e_1))}\circ\f|f| dm \leq m\big(
\f^{-1}\f_0(\cc(\e_1))\big)^{1/q}\|f\|_p
 $$
 and
 $$
\int\chi_{\f_0(\cc(\e_1))}\circ\f_0|f|dm \leq m\big(
\f_0^{-1}\f_0(\cc(\e_1))\big)^{1/q}\|f\|_p.
 $$
Let $\e_2>0$ be some small constant (to be determined later
in terms of $\e$). Using Proposition \ref{p.dimension} and
taking $\e_1$ and $\de$ sufficiently small we can make
 \begin{equation}
  \int_{\f_0(\cc(\e_1))}|\cl_\f
    f-\cl_{\f_0}f |dm\leq \e_2\var(f). \label{f1}
 \end{equation}

By a change of variables induced by $\f_0$ we deduce for
each $D_0\in\cad(\f_0)$ and $D\in\cad(\f)$ associated to
$D_0$
 \begin{eqnarray*}
 \lefteqn{\int_{\f_0(D_0)\cap
\f(D)}|\cl_\f f-\cl_{\f_0}f|dm=}\\
 &&
 \hspace{2cm}=
\int_{\widehat{D}_0}\left|\frac{f}{|\det
D\f|}\circ\f^{-1}\circ\f_0 -\frac{f}{|\det
D\f_0|}\right|\cdot|\det D\f_0|dm,
 \end{eqnarray*}
 and this
last expression is bounded from above by
 $$
 \int_{\widehat
D_0} \left(|f\circ\f^{-1}\circ\f_0-f|\cdot \frac{|\det
D\f_0|}{|\det D\f|\circ\f^{-1}\circ\f_0}+|f|\cdot
\left|\frac{|\det D\f_0|}{|\det
D\f|\circ\f^{-1}\circ\f_0}-1\right| \right)dm.
 $$
 Choosing
$\de>0$ sufficiently small, then $\|\f-\f_0\|_{C^1}<\de$
implies that
 $$
 \frac{|\det D\f_0|}{|\det
D\f|\circ(\f^{-1}\circ\f_0)}\leq 2 \quad\mbox{and}\quad
\left|\frac{|\det D\f_0|}{|\det
D\f|\circ(\f^{-1}\circ\f_0)}-1\right|\leq\e_2
 $$
 on
$S\setminus\cc(\e_1)$ (which contains $\widehat D_0$).
Hence
 $$
 \int_{\f_0(D_0)\cap \f(D)}|\cl_\f f-\cl_{\f_0}f|dm
\leq 2\int_{\widehat D_0} |f\circ\f^{-1}\circ \f_0-f|dm
+\e_2\int|f|dm.
 $$
 Thus, applying Lemma \ref{l.variation}
and choosing $\de>0$ sufficiently small we obtain
 \begin{equation}
 \int_{\f_0(D_0)\cap \f(D)}|\cl_\f
f-\cl_{\f_0}f|dm \leq \e_2\var(f).\label{f2}
 \end{equation}

Let us finally estimate the terms involved in (\ref{e3})
(the same method can be applied to obtain a similar
estimate for (\ref{e4})).
 \begin{eqnarray*}
\int_{\f_0(D_0)\setminus\f(D)}|\cl_{\f_0} f|dm &\leq &
\int_{\f_0(D_0)\setminus\f(D)}\frac{|f|} {|\det
D\f_0|}\circ\f_0^{-1}dm\\ &=&
\int_{D_0\setminus\f_0^{-1}(\f(D))}|f|dm\\ &\leq &
m(D\setminus \widehat D_0)^{1/q}\|f\|_p
 \end{eqnarray*}
 Using Proposition \ref{p.dimension} and  taking
 $\delta$ is sufficiently small, then
 \begin{equation}
\int_{\f_0(D_0)\setminus\f(D)}|\cl_{\f_0} f|dm+
\int_{\f(D)\setminus\f_0(D_0)}|\cl_{\f_0} f|dm \leq
\e_2\var(f).\label{f3}
 \end{equation}
  Putting estimates
(\ref{f1}), (\ref{f2}), (\ref{f3}) above together we obtain
 $$
 \int |\cl_\f f-\cl_{\f_0}f |dm \leq
\big(\e_2+3\#\cad(\f_0)\e_2\big)\var(f).
 $$
 So we only have
to take $\e_2$ in such a way that
$\e_2+3\#\cad(\f_0)\e_2<\e$.
\cqd

At this point we prove that conditions (U2') and (U2) are
equivalent if we assume (U1). First we prove that (U2')
implies (U2). Let $\e>0$ be some small number and take
$N\geq 1$ in such a way that
 $
 \|\sum_{j=N}^\infty\chi_{\{h_{\f_0}>j\}}\|_q<\e/3.
 $
 We have
 \begin{eqnarray*}
 \big\|\sum_{j=N}^\infty\chi_{\{h_\f>j\}}\big\|_q
 &= &
 \big\|h_\f-h_{\f_0}+
 h_{\f_0}-\sum_{j=0}^{N-1}\chi_{\{h_{\f_0}>j\}}
 +\sum_{j=0}^{N-1} \chi_{\{h_{\f_0}>j\}}-\sum_{j=0}^{N-1}
 \chi_{\{h_{\f}>j\}}\big\|_q\\
  &\leq &
 \|h_\f-h_{\f_0}\|_q+
 \big\|\sum_{j=N}^\infty\chi_{\{h_{\f_0}>j\}}\big\|_q
 + \sum_{j=0}^{N-1}\big\|\chi_{\{h_{\f_0}>j\}}-
 \chi_{\{h_{\f}>j\}}\big\|_q,
 \end{eqnarray*}
 and so, if we take $\de=\de(N,\e)>0$ sufficiently small
 then, under assumptions (U1) and (U2'),
 the first and third terms in the sum above can be made smaller
 than $\e/3$. This gives the conclusion of condition (U2).

 For the other implication, let $\e>0$ be some small number, and
 $N\geq 1$ and $\de=\de(N,\e)>0$ be taken in such a way that
 the conclusion of (U2) holds for $\e/3$. We have
 \begin{eqnarray*}
 \|h_\f-h_{\f_0}\|_q
 &=&
 \big\|h_\f
 -\sum_{j=0}^{N-1}\chi_{\{h_{\f}>j\}}
 +\sum_{j=0}^{N-1}\big(\chi_{\{h_{\f}>j\}}
 -\chi_{\{h_{\f_0}>j\}}\big)
 +\sum_{j=0}^{N-1}\chi_{\{h_{\f_0}>j\}}
 -h_{\f_0}\big\|_q \\
 &\leq&
 \big\|
 \sum_{j=N}^{\infty}\chi_{\{h_{\f}>j\}}\big\|_q
 +\sum_{j=0}^{N-1}\|\chi_{\{h_{\f}>j\}}
    -\chi_{\{h_{\f_0}>j\}}\|_q
 +
 \big\|\sum_{j=0}^{N-1}\chi_{\{h_{\f_0}>j\}}
 \big\|_q.
 \end{eqnarray*}
 By the choices of $N$ and $\de$,
 the first and third terms in the last sum above are
 smaller than $\e/3$. Since condition (U1)
 is verified, the second term can also be made
 smaller than $\e/3$ for $\de=\de(\e,N)>0$ small enough.

\cpr\label{p.cont*}
$\displaystyle{\frac{d\mu_{n_i}^*}{dm}}$ converges to
$\displaystyle{\frac{d\mu_{\infty}^*}{dm}}$ in the
$L^1$-norm.
 \fpr
  \dem Fixing some small $\e>0$, we are
going to prove that there is some $\de>0$ for which
$$
\left\|\frac{d\mu_{n_i}^*}{dm}-\frac{d\mu_{\infty}^*}{dm}
\right\|_1<\e \quad\mbox{whenever}\quad
\|\f_{n_i}-\f_0\|_{C^1}<\delta.
$$
We have
 \begin{equation}\label{e.formula}
  \mu^{\ast}_\infty
=\sum_{j=0}^{\infty}(\f_0^j)_*\left(\mu_\infty\mid
\{h_0>j\}\right)\qand \mu^{\ast}_{n_i}
=\sum_{j=0}^{\infty}(\f_{n_i}^j)_*\left(\mu_{n_i}\mid
\{h_{n_i}>j\}\right).
 \end{equation}
By (U2) we know that there is an integer $N\geq 1$ and
$\de=\de(\e,N)>0$ for which
 \begin{equation}\label{e.resto}
\|\f-\f_0\|_{C^1}<\de\Rightarrow
\big\|\sum_{j=N}^\infty\chi_{\{h_\f>j\}}\big\|_q<\frac{\e}{4K_1K_3}.
 \end{equation}
Now we take $i\geq 1$ sufficiently large in such a way that
$\|\f_{n_i}-\f_0\|<\de$. We split each one of the sums in
(\ref{e.formula}) into two sums and write
 \begin{equation}
\mu^{\ast}_\infty =\sum_{j=0}^{N}\nu_{\infty,j}+\eta_{\infty,N}\qand
\mu^{\ast}_{n_i} =\sum_{j=0}^{N}\nu_{n_i,j}+\eta_{n_i,N},
 \end{equation}
 where
 $$ \nu_{\infty,j}=\f_*^j\left(\mu_0\mid \{h_0>j\}\right),\quad
\eta_{\infty,N} =\sum_{j=N+1}^{\infty}\f_*^j\left(\mu_\infty\mid
 \{h_0>j\}\right)
 $$
and
 $$
 \nu_{n_i,j}=(\f_{n_i}^j)_*\left(\mu_{n_i}\mid
\{h_{n_i}>j\}\right),\quad
 \eta_{n_i,N}=\sum_{j=N+1}^{\infty}(\f_{n_i}^j)_*
 \left(\mu_{n_i}\mid\{h_{n_i}>j\}\right).
 $$
 We have
$$
 \eta_{\infty,N}(M)=\sum_{j=N}^\infty\mu_\infty(\{h_{0}>j\})
 =\sum_{j=N}^\infty\int\rho_{\infty}\chi_{\{h_{0}>j\}}dm
 \leq\|\rho_{\infty}\|_p\cdot\big\|
 \sum_{j=N}^\infty\chi_{\{h_{0}>j\}}\big\|_q,
 $$
and
 $$
 \eta_{n_i,N}(M)=\sum_{j=N}^\infty\mu_{n_i}(\{h_{n_i}>j\})
 =\sum_{j=N}^\infty\int\rho_{n_i}\chi_{\{h_{n_i}>j\}}dm
 \leq\|\rho_{n_i}\|_p\cdot\big\|
 \sum_{j=N}^\infty\chi_{\{h_{n_i}>j\}}\big\|_q
 $$
which together with Proposition \ref{p.dimension} and
  (\ref{e.resto}) yield
 \begin{equation}
 \left\|\frac{d\eta_{n_i,N}}{dm}-\frac{d\eta_{\infty,N}}{dm}\right\|_1
 \leq
 \eta_{n_i,N}(M)+\eta_{\infty,N}(M)<\e/2.
 \end{equation}
On the other hand, we have for $j=1,\dots,N$
 $$
 \left\|\frac{d\nu_{n_i,j}}{dm}-\frac{d\nu_{\infty,j}}{dm}\right\|_1
 =\big\|\cl_{\f_{n_i}^j}(\rho_{n_i}\chi_{\{h_{n_i}>j\}})-
 \cl_{\f_0^j} (\rho_{\infty}\chi_{\{h_{0}>j\}})\big\|_1
 $$
which is bounded from above by
 \[
 \begin{array}{ccc}
  \underbrace{\|\cl_{\f_{n_i}^j}(\rho_{n_i}\chi_{\{h_{n_i}>j\}})-
 \cl_{\f_{n_i}^j} (\rho_{\infty}\chi_{\{h_{0}>j\}})\|_1}&
 +&
 \underbrace{\|\cl_{\f_{n_i}^j}(\rho_{\infty}\chi_{\{h_{0}>j\}})-
 \cl_{\f_0^j} (\rho_{\infty}\chi_{\{h_{0}>j\}})\|_1}.\\
 A & & B
 \end{array}
 \]
 Here we also consider the transfer operator for the
 iterated maps $\f^j_{n_i}$ and  $\f^j_{0}$ defined
 analogously as for $\f$ in (\ref{e.transfer}).
We have
 \begin{eqnarray*}
 A
 & \leq&
 \|\rho_{n_i}\chi_{\{h_{n_i}>j\}}-
 \rho_{\infty}\chi_{\{h_{0}>j\}}\|_1\\
 & \leq&
 \|\rho_{n_i}\chi_{\{h_{n_i}>j\}}-
 \rho_{\infty}\chi_{\{h_{n_i}>j\}}\|_1
 +
 \|\rho_{\infty}\chi_{\{h_{n_i}>j\}}-
 \rho_{\infty}\chi_{\{h_{0}>j\}}\|_1\\
 &\leq&
 \|\rho_{n_i}-\rho_{\infty}\|_1
 +
 \|\rho_{\infty}(\chi_{\{h_{n_i}>j\}}-
 \chi_{\{h_{0}>j\}})\|_1
\end{eqnarray*}
 and
 $$
 \|\rho_{\infty}(\chi_{\{h_{n_i}>j\}}-\chi_{\{h_{0}>j\}})\|_1
 \leq
 \|\rho_{\infty}\|_p\|\chi_{\{h_{n_i}>j\}}-\chi_{\{h_{0}>j\}}\|_q.
 $$
 Taking into account (U1)
 we can make $A\leq \e/(4N)$ if $i$ is sufficiently large.
 Using Proposition \ref{l.converge} we can also make
 $
 B\leq \e/(4N),
 $
as long as we take $i$ large enough. This completes the
proof of the proposition.
\cqd

\cpr\label{p.igual}
 $\mu_\infty^*$ is a $\f_0$-invariant measure.
\fpr
 \begin{proof}
It follows from Proposition \ref{p.cont*} that
$(\mu_{n_i}^*)_i$ converges to $\mu_\infty^*$ in the weak*
topology. Hence, given any $f\colon M\rightarrow\RR$
continuous we have
 $$
\int fd\mu_{n_i}^*\rightarrow \int
fd\mu_{\infty}^*\quad\mbox{when}\quad i\rightarrow\infty.
 $$ On
the other hand, since $\mu_{n_i}^*$ is $\f_{n_i}$-invariant
we have
 $$
 \int fd\mu_{n_i}^*=\int
f\circ\f_{n_i}d\mu_{n_i}^*\quad\mbox{for every }i.
 $$
 It suffices to prove that
 \begin{equation}
 \int f\circ\f_{n_i}d\mu_{n_i}^* \rightarrow \int f\circ\f_0
d\mu_{\infty}^*\quad\mbox{when}\quad i\rightarrow\infty.
 \end{equation}
 We have
 \begin{eqnarray*}
 \lefteqn{
\big|\int f\circ\f_{n_i}d\mu_{n_i}^*-\int f\circ\f_0
    d\mu_{\infty}^*\big| \leq }\\
& &\hspace{2cm} \big|\int f\circ\f_{n_i}d\mu_{n_i}^*-\int
f\circ\f_0 d\mu_{n_i}^*\big|+\big|\int f\circ\f_0
d\mu_{n_i}^*-\int f\circ\f_0 d\mu_{\infty}^*\big|.
 \end{eqnarray*}
Since $f\circ\f_{n_i}-f\circ\f_0$ is uniformly close to zero when
$i$ is large (at least in the compact set $U$ that contains the supports of the measures), we have that the first term in the sum above
is close to zero for $i$ sufficiently large.
On the other hand, since $(\mu_{n_i}^*)_i$ converges to
$\mu_\infty^*$ in the weak* topology we also have that the second
term in the sum above
is close to zero if $i$ is taken large enough.
\end{proof}

It follows from this last result and the uniqueness of the
absolutely continuous $\f$-invariant measure that
necessarily $\mu_\infty^*=\mu^*_0$, thus proving that the measures
$\mu^*_{n_i}$ have densities converging in the $L^1$-norm
to the density of $\mu^*_0$. Moreover, the argument shows
that this  happens with the densities of any
convergent subsequence of $(\mu^*_n)_n$. This
completes the proof of Theorem~\ref{t.A}.


\section{Hyperbolic returns}\label{h.returns}
\label{s.hyperbolicreturns}

The main goal of this section is to introduce a notion of hyperbolic
returns, which allows us to improve some of the estimates in \cite{A}
and \cite{V} that are useful in the proofs of Theorem~\ref{t.B} and
Theorem~\ref{t.C}.
For the sake of clearness, we start by assuming that the map $\f $
has the special form
\begin{equation}
\f(\t,x)=(g(\t),f(\t,x)),\quad\mbox{with}\quad\partial_xf(\t,x)=0
\quad\mbox{if and only if}\quad x=0,
\label{assumption}
\end{equation}
and prove the conclusions of Theorems~\ref{t.A} and \ref{t.B}
for every $C^2$ map $\f $ satisfying
\begin{equation}
\|\f-\f_\al\|_{C^2}\leq\al\quad\mbox{on}\quad S^1\times I.\label{alfa}
\end{equation}
Later we explain how the conclusions extend to general case, using
the existence of a central invariant foliation in the same way as
in \cite{V} and \cite{A}.

Our estimates on the  derivative  depend in an
important way on the returns of orbits to the
neighborhood $S^1\times (-\sqrt{\al},\sqrt{\al}\,)$
of the critical set $\{x=0\}$. For this, we introduce
a partition $\cq $ of $I$ (modulo a zero Lebesgue measure set)
into the following intervals:
$$
I_r=(\sqrt{\al} e^{-r},\sqrt{\al} e^{-(r-1)}\,)\quad \mbox{for}\quad r\geq 1,
$$
$$
I_r=-I_{-r}\quad \mbox{for}\quad r\leq -1,
$$
$$
I_{0^+}=\left(I\setminus[-\sqrt{\al},\sqrt{\al}\,]\right)\cap\RR^+
\quad \mbox{ and }\quad
I_{0^-}=\left(I\setminus[-\sqrt{\al},\sqrt{\al}\,]\right)\cap\RR^-.
$$
This partition induces in a natural way
analogous ones at each fiber of the type
$\{\t\}\times I$. For the sake simplicity in the notation
no specification will be made in which fiber
they are on, since this will be always clear in all
our settings.

In what follows we assume that $\al >0$ is a sufficiently
small number independent of any other constant involved in
the arguments. Furthermore, for each new constant
appearing  we will always specify whether it depends on
$\al$ or not. Given  $(\t,x) \in S^1\times I$ and $j\geq 0$
we define $(\t_j,x_j) =\f^j(\t,x)$. Following \cite{V}, for
the next lemma we take $\eta$ a positive constant smaller
than $1/4$ depending only on the quadratic map $q$.

\cle \label{N}
There are constants $C_1>1$ such that for every small $\al $
we have an integer $N(\al)$ satisfying:
\begin{enumerate}
\item If $|x|<2\sqrt{\al} $, then
$\prod^{N(\al)-1}_{j=0}|\partial_xf(\t_j,x_j)|\geq|x|\al^{-1+\eta}$.
\item If $|x|<2\sqrt{\al} $, then $|x_j|>\sqrt{\al} $ for every
$j=1,\dots,N(\al)$.
\item $C_1^{-1}\log(1/\al)\leq N(\al)\leq C_1\log(1/\al)$.
\end{enumerate}
\fle
\dem See \cite[Lemma 2.4]{V} and \cite[Lemma 2.1]{A}.\cqd

\cle\label{lema1}
There are $\tau>1$, $C_2>0$ and $\de>0$ such that for $(\t,x)\in S^1\times I$
and $k\geq1$ the following holds:
\begin{enumerate}
\item If $|x_0|,\dots,|x_{k-1}|\geq\sqrt{\al} $, then
    $\prod^{k-1}_{j=0}|\partial_xf(\t_j,x_j)|\geq C_2\sqrt{\al}\tau^k$.
\item If $|x_0|,\dots,|x_{k-1}|\geq\sqrt{\al} $ and $|x_k|<\de $,  then
    $\prod^{k-1}_{j=0}|\partial_xf(\t_j,x_j)|\geq C_2\tau^k$.
\end{enumerate}
\fle \dem See \cite[Lemma 2.5]{V}. \cqd

For each integer $j\geq 0$ we define
\begin{equation}
\label{e.defr}
r_j(\t,x)=\left\{
\begin{array}{cl}
|r| & \mbox{ if }\f^j (\t,x)\in I_r\mbox{ with }|r|\geq 1;\\
0 & \mbox{ if }\f^j (\t,x)\in I\setminus [-\sqrt{\al},\sqrt{\al}\,].
\end{array}\right.
\end{equation}
We say that $\nu\geq 0 $ is a {\it return} for $(\t,x)$
if $r_\nu(\t,x) \geq 1$.
Let $n$ be some positive integer and
$0\leq\nu_1\leq\cdots\leq\nu_s\leq n$
the returns of $(\t,x)$ from $0$ to  $n$.
It follows from Lemma \ref{N} that for each
$1\leq i\leq s$
$$\prod_{j=\nu_i}^{\nu_i+N-1}|\partial_xf(\t_j,x_j)|
     \geq e^{-r_{\nu_i}(\t,x)}\al^{-1/2+\eta },$$
and from Lemma \ref{lema1}
$$\prod_{j=0}^{\nu_{1}-1}|\partial_xf(\t_j,x_j)|
     \geq C_2\tau^{\nu_{1}}\quad\mbox{and}\quad
\prod_{j=\nu_i+N}^{\nu_{i+1}-1}|\partial_xf(\t_j,x_j)|
     \geq C_2\tau^{\nu_{i+1}-\nu_i-N}.$$
For the last piece of orbit (if it exists) we use again Lemma \ref{lema1}
and obtain
\begin{equation}\label{e.ultima}
\prod_{j=\nu_s}^{n-1}|\partial_xf(\t_j,x_j)|
     \geq |\partial_x f(\t_{\nu_s},x_{\nu_s})|C_2\sqrt{\al}\tau^{n-\nu_{s}}.
\end{equation}
Considering
\begin{equation}
\label{e.defG}
G_n(\t,x)=\left\{1\leq \nu_i\leq n-1:r_{\nu_i}(\t,x)
\geq \left(\frac{1}{2}-2\eta\right)\log\frac{1}{\al}\right\},
\end{equation}
the estimates above yield (see \cite[Section 2]{A})
\begin{equation}\label{e.ineq1}
\prod_{j=0}^{n-1}|\partial_xf(\t_j,x_j)|
\geq \exp\left(4cn-\sum_{j\in G_n(\t,x)}r_{j}(\t,x)-\log C_2
 -\frac{3}{2}\log\frac{1}{\al}\right)
\end{equation}
for every $n\geq 1$ and $\al $ sufficiently small, where
$c>0$ is some constant depending only on the map $q$.
The term $(3/2)\log(1/\al)$ appears if
a last piece of orbit
has to be considered whenever
$n$ is not a return for $(\t,x)$
(estimate (\ref{e.ultima}) above). Hence, if $n$ is a return
for $(\t,x)$ we can improve estimate (\ref{e.ineq1})
and deduce that
\begin{equation}\label{e.ineq2}
\prod_{j=0}^{n-1}|\partial_xf(\t_j,x_j)|
\geq C_2^{-1}\exp\left(4cn-\sum_{j\in G_n(\t,x)}r_{j}(\t,x)
 \right)
\end{equation}

A key fact underlying our construction is that the exponent
on the right hand side is positive, except for a set of
initial points $(\t,x)$ whose measure decreases very
rapidly with $n$. More precisely, let us define
\begin{equation} \label{e.defE}
 E_n=\big\{(\t,x)\in S^1\times
I: \sum_{j\in G_n(\t,x)}r_{j}(\t,x) > 2n\big\}.
 \end{equation}
Then, cf. (16) and (17) in \cite[Section 2.4]{V}, there are
constants $C,\gamma>0$ such that
\begin{equation}\label{e.med}
m\big(E_n\big) \leq Ce^{-\ga\sqrt n}
\end{equation}
for every sufficiently
large $n$, only depending on $\al$. Note that
(\ref{e.ineq2}) gives
 $$
 \|D\f^n(\theta,x)(0,1)\| \ge e^{cn}
\quad\text{for $(\t,x)\in (S^1\times I) \setminus E_n$}\,
 $$
and $n$ sufficiently large.
\medskip

One of the basic ingredients in the proof of the existence of the SRB
measures is the notion of hyperbolic times.
Following \cite{A} we fix $0<\e<c/2$ and say that $n\geq 1$ is a {\em
hyperbolic time} for $(\t,x)\in S^1\times I$ if
$$
\sum_{\mbox{\scriptsize $
\begin{array}{c}
i\in G_n(\t,x)\\
k\leq i<n
\end{array}$}
}
r_i(\t,x)<(c+\e)(n-k)\quad\mbox{for every}\quad 0\leq k< n.
$$

We say that $n\geq 1$ is a {\em hyperbolic return}
for $(\t,x)\in S^1\times I$ if $n$ is both a
hyperbolic time and a return for $(\t,x)$. It follows from
\cite[Proposition 2.5]{A} (see also \cite[Remark 2.6]{A})
that Lebesgue almost every point in $S^1\times I$ has
infinitely many hyperbolic times. This in particular
implies that Lebesgue almost every point in $S^1\times I$
also has infinitely many hyperbolic returns. Indeed, if $n$
is a hyperbolic time for $(\t,x)$ and $l>n$ is the next
return for $(\t,x)$ after $n$, then since $r_{j}(\t,x)=0$
for $j=n+1,\cdots,l-1$, it easily follows that $l$ is a
hyperbolic return for $(\t,x)$.

Fixing an integer $p\geq 1$ (whose value will be made
precise below in terms of the expansion rates of
the maps $\hat g$ and $\hat f$), let $H$ be the set of points
that has at least one hyperbolic time greater or equal to
$p$. Decompose $H=\cup_{n\geq p}H_n$, where each $H_n$
is the set of points whose first hyperbolic time greater or
equal to $p$ is $n$. It follows from \cite[Proposition
2.5]{A} that there is a positive integer $n_0=n_0(p,\e)\geq p$
such that
\begin{equation}\label{e.inc}
(S^1\times I)\setminus (H_p\cup\cdots\cup H_n)\subset E_n
\quad\mbox{for every}\quad n\geq n_0\,.
\end{equation}

Now we briefly describe how in \cite[Section 3]{A} is
defined a partition $\car$ into rectangles of $S^1\times I$
(modulo a zero Lebesgue measure set). For this, we consider
the partition $\cq$ of $I$ described above, and introduce a
sequence of Markov partitions of $S^1$: assume that
$S^1=\RR/\ZZ$ has the orientation induced by the usual
order in $\RR$ and let $ {\t}_0$ be the fixed point of $g$
close to $\t=0$. We define Markov partitions $\cp_n$,
$n\geq 1$, of $S^1$ in the following way:
\begin{itemize}
\item $\cp_1=\{[ {\t}_{j-1}, {\t}_{j}):1\leq j\leq d\}$,
where $ {\t}_0, {\t}_1,\dots, {\t}_d= {\t}_0$
are the pre-images of $ {\t}_0$ under $g$ (ordered according
to the orientation of $S^1$).
\item $\cp_{n+1}=\{\mbox{ connected components of }g^{-1}(\o):
\o\in\cp_n\}$ for each $n\geq 1$.
\end{itemize}
The partition $\car $ is obtained by successive divisions
of the rectangles in an initial partition $\cp_p\times\cq$
of $S^1\times I$, for some fixed large integer $p$,
according to the itineraries of points through the
horizontal strips $S^1\times I_*$ with $I_*\in\cq$ and
their hyperbolic times. This partition may be written as a
union $\car=\cup_{n\geq p}\car_n $ with the sets $\car_n$
defined inductively and satisfying
\begin{equation}\label{cobre}
H_n\subset\bigcup_{R\in\car_n}R\quad\mbox{and}\quad
R\cap H_n\neq\emptyset\quad\mbox{for every}\quad R\in\car_n.
\end{equation}
For $n\geq p$ rectangles in $\car_n $ always have the form
$\o\times J$, with $\o $ belonging to  $\cp_n$ and $J$
a subinterval of $I_*$ for some $I_*\in\cq $.

\cpr \label{p.dist0}
There is some constant $\De>1$ such that for every $n\geq p$,
$R\in \car_n$ and $(\t,x),(\si,y)\in R$ we have
$$
\frac{1}{\De}\leq \left|\frac{J(\t,x)}{J(\si,y)}\right|\leq \De,
$$
where $J$ is the Jacobian of $\f^n\mid R$.
\fpr
\dem
Fix some
$R\in\car_n$ with $n\geq p$ and let $\F=\f^{n}\mid R$. We have
$$ \left|\frac{J(\t,x)}{J(\si,y)}\right|=
\exp\left(\log\left|\big(J\circ\F^{-1}\big)\left(
\F(\t,x)\right)\right|-
          \log\left|\big(J\circ\F^{-1}\big)\left(
\F(\si,y)\right)\right|\right)
$$
and
$$
\left|\log\left|\big(J\circ\F^{-1}\big)\left(
\F(\t,x)\right)\right|-
          \log\left|\big(J\circ\F^{-1}\big)\left(
\F(\si,y)\right)\right|\right|\leq
\left\|D\left(\log\left|J\circ\F^{-1}\right|\right)
(\tau,z)\right\|\cdot C,$$
for some $(\tau,z)\in \F(S)$ and $C>0$ depending only
on the diameter of $S^1\times I$.
Now, since
$$\left\|D\left(\log\left| J\circ\F^{-1}\right|\right)
(\tau,z)\right\|=
\frac{\left\|D\left(J\circ\F^{-1}\right)
(\tau,z)\right\|}{\left|(J\circ\F^{-1})(\tau,z)\right|},$$
the result follows from \cite[Proposition 4.2]{A}.
\cqd

Now we are going to prove that the Lebesgue measure of the
set of points that have no hyperbolic returns smaller than
some large integer $n$ decays at least sub-exponentially
fast with $n$. Similarly to what we have done for
hyperbolic times, let $p\geq 1$ be some fixed large
integer, and define $H^\ast$ the set of points that has at
least one hyperbolic return greater or equal to $p$. We
decompose $H^\ast=\cup_{n\geq p}H_n^\ast$, where each
$H_n^\ast$ is the set of points whose first hyperbolic
return greater or equal to $p$ is $n$.

\cpr \label{p.cauda}
There is a positive integer $n_1=n_1(p,\e)\geq p$
and constants $C_0,\gamma_0>0$ such that for each $n\geq
n_1$ $$ m\big((S^1\times I)\setminus (H^*_p\cup\cdots\cup
H^*_{n})\big)\leq C_0e^{-\gamma_0\sqrt n}. $$
\fpr
\dem
Take $n\geq \max\{2p,n_0\}$ and let $l=[n/2]$. The set of points
$(\t,x)\in H_l$ for which there is some $1\leq k\leq l$
that is a return for $\f^l(\t,x)$ is contained in
$H^*_p\cup\cdots\cup H^*_{n}$. Hence,
defining
 $$ B_l=\bigcup_{k=p}^l\{(\t,x)\in H_k\colon
\f^k(\t,x) \mbox{ has no returns from time 1 to $l$} \}
 $$
we have
$$
(H_p\cup\cdots\cup H_l)\cap \big((S^1\times I)\setminus B_l\big)
\subset (H_p^*\cup\cdots\cup H_n^*)
$$
and so
\[
m\big((S^1\times I)\setminus (H^*_p\cup\cdots\cup
H^*_{n})\big) \leq m\big((S^1\times I)\setminus
(H_p\cup\cdots\cup H_l)\big)+m(B_l),
\]
Taking into account estimates (\ref{e.inc}) and
(\ref{e.med}) above, it
suffices to study the decay of $m(B_l)$ with $n$. We define
for each $k\geq p$ and $R_k\in\car_k$
 $$
  R_k(l)=\{(\t,x)\in R_k\colon
\f^k(\t,x) \mbox{ has no returns from time 1 to $l$} \}. $$
Using (\ref{cobre}) we obtain
\begin{equation}\label{zero}
m(B_l)\leq \sum_{k=p}^l\sum_{R_k\in
\car_k}m\big(R_k(l)\big).
\end{equation}
Fixing some $R_k\in\car_k$ and $(\t_0,x_0)\in R_k$  we
deduce from Proposition \ref{p.dist0}
\begin{eqnarray*}
m\left(\f^{k}\big(R_k(l)\big)\right) &=&
\int_{R_k(l)}\left|J(\t,x)\right|dm(\t,x)\\ &\geq
&\frac{1}{\De}\int_{R_k(l)}\left|J(\t_0,x_0)\right|dm(\t,x)\\
&\geq
&\frac{1}{\De}\left|J(\t_0,x_0)\right|m\big(R_k(l)\big).
\end{eqnarray*}
Similarly we prove that
 $$
 m\left(\f^{k}(R_k)\right) \leq
{\De}\left|J(\t_0,x_0)\right|m(R_k).
 $$ Hence
\begin{equation}\label{e.um}
\frac{m\big(R_k(l)\big)}{m(R_k)}\leq\De^2
\frac{m\left(\f^k\big(R_k(l)\big)\right)}{m\left(\f^k(R_k)\right)}.
\end{equation}
It follows from the definition of $R_k(l)$ that the
iterates of points in $\f^k\big(R_k(l)\big)$ do not hit the
critical region $S^1\times [-\sqrt\al,\sqrt\al]$ from time
$1$ to $l$. From Lemma \ref{lema1} we deduce that there is
some constant $C>0$ for which
\begin{equation}\label{e.dois}
m\left(\f^k\big(R_k(l)\big)\right)\leq C\tau^{-l}.
\end{equation}
On the other hand, it follows from \cite[Proposition 3.8]{A}
that there is some absolute constant $\de>0$ for which
\begin{equation}\label{e.tres}
m\left(\f^k\big(R_k\big)\right)\geq \de.
\end{equation}
From (\ref{e.um}), (\ref{e.dois}) and (\ref{e.tres}) we
obtain
 $$ m\big(R_k(l)\big)\leq \frac{\De^2
C}{\de}\tau^{-l} m(R_k)
 $$
 which together with (\ref{zero})
gives
 $$ m(B_l)\leq \sum_{k=p}^l\sum_{R_k\in \car_k}
\frac{\De^2 C}{\de}\tau^{-l}m(R_k) \leq (l-p)\frac{\De^2
C}{\de}\tau^{-l}\leq n\frac{\De^2 C}{\de}\tau^{-n/2}.
 $$
\cqd

\cre \label{r.cgama}
It follows from the proof of Proposition \ref{p.cauda}
that the constants $C_0$ and $\gamma_0$ only depend
on the constants $C$, $\gamma$ and absolute constants
associated to the quadratic map $q$.
Moreover, the integer $n_1$ only depends on the previous
constants and the integer
$p\geq 1$. At the end of this section we will see that
$p$ may be chosen independent of the map $\f\in\cn$.
\fre

Hyperbolic times play a crucial role in \cite[Proposition 3.8]{A}
in order to obtain that the images $\f^n(R)$ of rectangles
$R\in\car_n $ have sizes uniformly bounded away from zero.
However, the uniform constant that bounds such sizes from
below depends on $\al$, and that is still inconvenient for
proving the mixing and ergodic properties.
 To bypass this difficulty we are going to consider
hyperbolic returns in the place of hyperbolic times. We
proceed as in \cite[Section 3]{A} and define a partition
$\car$ into rectangles of $S^1\times I$ (modulo a zero
Lebesgue measure set) exactly in the same way
with the sets $H^*_n$ playing the role of the sets $H_n$.
In particular, this partition may also be written as a
union $\car=\cup_{n\geq p}\car_n $ with the sets $\car_n$
defined inductively and satisfying $$
H_n^*\subset\bigcup_{R\in\car_n}R\quad\mbox{and}\quad R\cap
H_n^*\neq\emptyset\quad\mbox{for every}\quad R\in\car_n. $$
Furthermore, for each $n\geq p$, rectangles in $\car_n $
also have the form $\o\times J$, with $\o $ belonging to
$\cp_n$ and $J$ a subinterval of $I_*$ for some $I_*\in\cq
$. We define a map $h:\car\rightarrow \ZZ^+$, by putting
$h(R)=n\geq p$ for each $R\in\car_n$.

\cle\label{l.exp}
Let $(\t,x)\in R$ for some $R\in\car$.
Then for every $j=0,\cdots, h(R)-1$ we have
$$
\prod_{i=j}^{h(R)-1}|\partial_x f(\t_i,x_i)|\geq
C_2^{-1}\exp \big((2c-\e)(h(R)-j)\big).
$$
\fle
\dem The same proof of \cite[Lemma 3.7]{A} with
the improved estimate (\ref{e.ineq2}) in the place
of (\ref{e.ineq1}). \cqd

It follows from assumption (\ref{assumption}) that
for each $n\geq 1$ there is a map
$F_n\colon S^1\times I\rightarrow I$
such that $\f^n(\t,x)=(g^n(\t) , F_n(\t,x) )$
for every $(\t,x)\in S^1\times I$.
Let $(\t,x)$ belong to $R\in \car$ and $h=h(R)$. Then
\[
D\f^{h}(\t,x)=\left(
\begin{array}{cc}
       \partial_\t g^{h}(\t ) & 0 \\
\partial_{\t}F_{h} (\t,x)  & \partial_x F_{h}(\t,x)
\end{array}
\right),
\]
and so
\begin{eqnarray*}
(D\f^{h}(\t,x))^{-1}&=&\frac{1}{\partial_\t g^{h}(\t )\partial_x F_{h}(\t,x) }
\left(
\begin{array}{cc}
       \partial_x F_{h}(\t,x) & 0 \\
      -\partial_{\t}F_{h} (\t,x)  & \partial_\t g^{h}(\t )
\end{array}
\right)\\
& & \\
&=&\left(
\begin{array}{cc}
  (\partial_\t g^{h}(\t ))^{-1} & 0 \\
 \!-{\partial_{\t}F_{h} (\t,x)}\big({\partial_\t g^{h}(\t )\partial_x
F_{h}(\t,x)}\big)^{-1}
                       &({\partial_x F_{h}(\t,x)})^{-1}
      \end{array}
\right).
\end{eqnarray*}
It follows from \cite[Lemma 4.1]{A} that there is some
constant $C_3>0$ such that for every $(\t,x)\in S^1\times I$
we have
$|\partial_{\t}F_{h}(\t,x)|\leq C_3|\partial_{\t}g^{h}(\t)|$. Then
$$\|D\f^{-h}(\f^{h}(\t,x))\|\leq
\max\left\{|\partial_\t g^{h}(\t )|^{-1}+C_3|
\partial_x F_{h}(\t,x)|^{-1},
|\partial_x F_{h}(\t,x)|^{-1}\right\}.$$
We have
$$|\partial_{\t} g^{h}(\t)|^{-1}\leq (d-\al)^{-h}$$
and from Lemma \ref{l.exp}
$$|\partial_x F_{h}(\t,x)|^{-1}\leq C_2\exp(-(2c-\e)h).$$
Hence,
\begin{equation}
\|D\f^{-h}(\f^{h}(\t,x))\|\leq (d-\al)^{-h}+ (1+C_3)C_2
\exp(\!-(2c-\e)h),\label{norma}
\end{equation}
At this point we can specify the choice of the integer
$p$: we take $p\geq 1$ large enough in such that the induced
map $\F$ associated to $\f$ is an expanding map
in the sense of the definition given in
Subsection \ref{s.statement} (recall that $h\geq p$).
Note that this choice of $p$ only depends on the expansion rates of
the maps $\hat g$ and $\hat f$, thus $p$ may be taken independent of the map $\f\in\cn$.


\section{Uniformity conditions}


An important feature of this construction, cf. \cite[Section 2.5]{V},
is that it remains valid for any map $\psi$ close enough to $\f$,
with uniform bounds on the measure of the exceptional sets $E_n(\psi)$:
$$
m(E_n(\psi)) \le Ce^{-\gamma\sqrt{n}} \quad\text{for every } n\ge 1
$$
where $C$ and $\gamma$ may be taken uniform (that is, constant) in
a whole $C^3$ neighbourhood of $\f$.
Let us explain this last point, since it is not explicitly addressed
in the previous papers.
One consequence is that Proposition \ref{p.cauda} holds in the whole
open set $\cn$, with uniform constants $C_0$ and $\gamma_0$ (recall
Remark \ref{r.cgama}).

As explained in \cite[Section 2.5]{V}, it follows from the methods of
\cite{HPS} that any map $\psi$ sufficiently close to $\f$ admits a
unique invariant central foliation $\cf^c$ of $S^1\times I$ by smooth
curves uniformly close to vertical segments.
This is because the vertical foliation is invariant and normally
expanding for the map $\f$.
In addition, the space of leaves of $\cf^c$ is homeomorphic to a circle,
and the map induced by $\psi$ in it topologically conjugate to $\hat{g}$.
The previous analysis can then be carried out in terms of the
expansion of $\psi$ along this central foliation $\cf^c$.
More precisely, $|\partial_x f(\t,x)|$ is replaced by
$$
|\partial_c f(\t,x)| \equiv |D\psi(\t,x)v_c(\t,x)|,
$$
where $v_c(\t,x)$ represents a norm $1$
vector tangent to the foliation at $(\t,x)$.
The previous observations imply that $v_c$ is uniformly close to
$(0,1)$ if $\psi$ is close to $\f$.
Moreover, cf. \cite[Section 2.5]{V}, it is no restriction to
suppose $|\partial_c f(\t,0)|\equiv 0$
(incidentally, this is the only place where we need our maps
to be $C^3$), so that $\partial_c f(\t,x)\approx |x|$, as in
the unperturbed case; recall (\ref{assumption}).
Defining $r_j(\t,x)$ and $E_n=E_n(\psi)$ in the same way as
before, cf. (\ref{e.defr}), we obtain an analog of (\ref{e.ineq2}).
\begin{equation*}
\|D\psi^n(\theta,x)v_c(\theta,x)\|=
\prod_{j=0}^{n-1}|\partial_c f(\t_j,x_j)|
\geq C_2^{-1}\exp\left(4cn-\sum_{j\in G_n(\t,x)}r_{j}(\t,x) \right),
\end{equation*}
for every $(\t,x)$.
We define $E_n(\psi)$ in the same way as $E_n=E_n(\f)$,
recall (\ref{e.defE}), and then
\begin{equation*}
\|D\psi^n(\theta,x)v_c(\theta,x)\| \ge e^{cn}
\quad\text{for all } (\t,x)\in (S^1\times I) \setminus E_n.
\end{equation*}

The arguments in \cite[Section 2.4]{V} apply with
$|\partial_c f|$ in the place of $|\partial_x f|$, proving
that the Lebesgue measure of $E_n(\psi)$ satisfies the
bound in (\ref{e.med}). The constants $C$ and $\gamma$
produced by these arguments depend only on $\alpha$, which
is fixed, and on estimates obtained in the previous
sections of that paper. So, to see that these constants are
indeed uniform in a neighbourhood of $\f$, it suffices to
check that the same is true for those preparatory
estimates. This is clear in the case of the results of
Section~2.1 (Lemmas 2.1 and 2.2, and Corollary 2.3),
because they only involve one iterate of the map. Let us
point out that the definition of admissible curve for
$\psi$ is just the same as for the unperturbed map $\f$. A
continuity argument can be applied also to Section 2.2, but
it is more subtle. The key observation is that, although
the statements of Lemmas 2.4 and 2.5 involve an unbounded
number of iterates, their proofs are based on analyzing
bounded stretches of orbits. Finally, the results in
Section 2.4 (Lemmas 2.6 and 2.7), involve not more than
$M\approx\log(1/\alpha)$ iterates. So, once more by
continuity, their estimates remain valid in a neighbourhood
of $\varphi$. We have concluded the observation that the
bound (\ref{e.med}) on the Lebesgue measure of the
exceptional set $E_n$ holds uniformly in a neighbourhood of
the map.



\medskip
Now we are able to show that conditions (U1)-(U3)
are satisfied by every element of $\cn$, as long as
we take the open set $\cn$ sufficiently small.

\begin{itemize} \item[(U1)] The construction of the
partition that leads to the map $h_\f$ is based on the
itineraries of points through the horizontal strips
$S^1\times I_*$ with $I_*\in \cq$, according to the
expanding behaviour of the iterates of $\f$ at hyperbolic
returns. Since these hyperbolic returns depend only on a
finite number of iterates of the map $\f$, by continuity,
we can perform the construction of the partition in such a
way that for some fixed integer $N$ the Lebesgue measure of
$\{h_\f=j\}$ varies continuously with the map $\f$ for
$j\leq N$.

\item[(U2)] We have for every $\f\in\cn$ and any fixed
large integer $N\geq 1$ $$ \big\|\sum_{j\geq
N}\chi_{\{h_\f>j\}}\big\|_q^q =\sum_{j\geq
N}j^qm\big(\{h_\f=N+j\}\big) \leq\sum_{j\geq
N}j^qm\big(\{h_\f\geq N+j\}\big) $$ Taking into account
Proposition \ref{p.cauda} we deduce that $$
\big\|\sum_{j\geq N}\chi_{\{h_\f>j\}}\big\|_q^q
\leq\sum_{j\geq N}j^qC_0e^{-\ga_0\sqrt{N+j}} $$ which can
be made uniformly small if $N$ is taken sufficiently large.

 \item[(U3)] The constant $K$
that bounds the distortion is given by \cite[Proposition
4.2]{A}, which may be taken uniform in the whole $\cn$.

The constant $\si$ is given by (\ref{norma}), which may
be taken uniformly smaller than one for every $\f\in \cn$.

It follows from \cite[Corollary 3.3]{A} that $\beta$ is
uniformly bounded away from zero, as long as $\al$ and the
open set $\cn$ are taken small enough.

Finally, the proof \cite[Proposition 3.8]{A} shows that
$\rho$ may be taken bounded from below by a constant only
depending on $\alpha$ (see also Step 1 in the proof of
Proposition \ref{p.whole} below).

\end{itemize}


\begin{Remark}
The following comments are meant to help clarify the presentation of
\cite[Lemma 2.6]{V}, they are not used elsewhere in the present work.
We refer the reader to \cite{V} for the setting and notations.
The conclusion of the lemma is contained in Corollary 2.3 of \cite{V},
when $r$ is large enough so that $|J(r-2)|\ll \sqrt{\alpha}$.
In particular, it is enough to prove the lemma for values of $r$
smaller than $(1/2+2\eta)\log(1/\alpha)$
(and larger than $(1/2-2\eta)\log(1/\alpha)$, cf. statement of the lemma).
By definition, the function $k(r)$ defined in page 73 of \cite{V}
can not exceed $M\approx\log(1/\alpha)$.
So, the arguments at the end of page 73 actually prove that either
$k(r) \ge \text{const\,} r$, or $k(r)\approx M$.
However, under the above restriction on $r$, the latter possibility
also implies $k(r)\ge \text{const\,} r$.
In this way, the conclusion of the lemma follows in all the cases.
\end{Remark}


\section{Topological mixing}

In this section we prove that the maps in  $\cn$ are
topologically mixing. For that, let us start by giving a
good description of the attractor of a map $\f\in\cn$
inside the forward invariant region $S^1\times I$.  We
claim that the attractor of $\f$ inside this invariant
region, defined as the intersection
 $$ \Lambda
=\bigcap_{n\geq 0}\f^n(S^1\times I)
 $$
 of all forward
images of $S^1\times I$, is just $\Lambda=\f^2(S^1\times
I)$, as long as the interval $I$ is properly chosen.
Indeed, for the one-dimensional map $q$ we may take
$I\subset (-2,2)$ in such a way that $q(I)$ is contained in
the interior of $I$ and
 $$
 \bigcap_{n\geq 0}q^n(I)= q^2(I).
 $$
 For each $\t\in S^1$ the map $\f\mid(\{\t\}\times I)$
may be thought of as a one-dimensional map from $I$ into
itself, close to $q$. Then, our claim follows by
continuity.

Before we go into the main proposition of this section, let
us remark that inequality (\ref{norma}) shows in particular
that the diameter of the partition $\car$ of $S^1\times I$,
defined as
 $$
 \mbox{diam}(\car)=\sup\,\{\mbox{diam}(R):R\in\car\},
 $$
 is small when $p$ is large. Thus, taking arbitrarily
large integers $p$ we may define  a sequence of partitions
$(\cs_n)_n$ in $S^1\times I$  in such a way that
 $$
 \lim_{n\rightarrow+\infty}\mbox{diam}\,(\cs_n)=0.
 $$
Moreover, for each $n\geq 1$ we may define a map
 $h_n:\cs_n\rightarrow \ZZ^+$
in the same way as we did for $h:\car\rightarrow \ZZ^+$.

\cpr \label{p.whole}
There is  integer $M=M(\al)$ such that
for every $n\geq 1$ and $\o\times J\in \cs_n$,
$$
\left|\f^{h_n(\o\times J)+M}(\o\times J)\right|=\La.
$$
\fpr
\dem The proof of this proposition will be made in
four steps. In the first one we prove that the image of
$\o\times J\in \cs_n$ by $\f^{h_n(\o\times J)}$ has height
bounded from below by a constant of order $\al^{1-2\eta}$.
In the second step we prove that a vertical segment
of order $\al^{1-2\eta}$ becomes, after
a finite number of iterates, an interval
with height bounded from below by a constant of order $\sqrt{\al}$.
In the third step we show that iterating vertical segments of
order $\sqrt{\al}$ they become segments with length bounded
from below by a constant not depending on $\al$. In the
final step we make use of the properties of the quadratic
map $q$ to obtain the
result.

\paragraph{Step 1.}
{\em There is a constant $\De_1>0$ such that for
every $n\geq 1$, $\o\times J\in \cs_n$ and $\t\in\o$},
$$\left|\f^{h_n(\o\times J)}(\{\t\}\times J)\right|
\geq \De_1\al^{1-2\eta}.$$

\smallskip
\noindent This follows from \cite[Proposition 3.8]{A}. It
is easy to check that the estimate in Lemma \ref{l.exp}
above in the place of the estimate of \cite[Lemma 3.7]{A}
is enough to yield this dependence on $\al$.

\paragraph{Step 2.}
{\em There is a constant $\De_2>0$ and an integer $M_1=M_1(\al)$
such that if $J\subset I$ is an interval with
$|J|\geq \De_1\al^{1-2\eta}$, then for every $\t\in S^1$}
$$
\left|\f^{M_1}(\{\t\}\times J)\right|\geq \De_2\sqrt\al.
$$

\smallskip
\noindent We start by remarking that we may assume that $J$
intersects $(-\sqrt\al,\sqrt\al)$. (If this is not the
case, we take $R\geq 1$  the first integer for which
$\f^{R}(\{\t\}\times J) $ intersects
$(-\sqrt\al,\sqrt\al)$. It follows from Lemma \ref{lema1}
that $R\leq C\log(1/\al)$ for some $C>0$, and so we may
start with $\f^{R}(\{\t\}\times J) $). Take $J_1$ a
subinterval of $ J$ such that $$ J_1\subset
(-2\sqrt\al,2\sqrt\al),\quad J_1\cap
\left(-\frac{\De_1}{4}\al^{1-2\eta},
\frac{\De_1}{4}\al^{1-2\eta}\right)=\emptyset \quad \mbox{
and }\quad |J_1|\geq \frac{\De_1}{4}\al^{1-2\eta}. $$ It
follows  from Lemma \ref{N} that  for $N=N(\al)$
\begin{equation}\label{e.step2}
\left|\f^N(\{\t\}\times J_1)\right|\geq \frac{\De_1}{4}\al^{1-2\eta}\;
\al^{-1+\eta}|J_1|
\geq \frac{\De_1}{4} \al^{1-3\eta}.
\end{equation}
Let $R\geq 1$ be the first integer for which
$\f^{N+R}(\{\t\}\times J_1) $ intersects
$(-\sqrt\al,\sqrt\al)$. It follows from Lemma \ref{lema1}
that $$ \left|\f^{N+R}(\{\t\}\times J_1)\right|\geq
\frac{C_2\De_1}{4} \al^{1-3\eta}. $$ Now we proceed
inductively and prove that for each $l\geq 1$ there is an
interval $J_l\subset J$ and a sequence of integers
$3=k_1<k_2<\cdots<k_l$ for which
$$
\left|\f^{lN+(l-1)R}(\{\t\}\times J_l)\right|\geq
    \frac{C_2^{l-1}\De_1}{4^l} \al^{1-k_l\eta}.
$$
We stop when $1-k_l\eta\leq 1/2$, and take
$\De_2={C_2^{l-1}\De_1}/{4^l} $ and $M_1=lN+(l-1)R$
(note that $l$ only depends on $\eta$ which does not depend
on $\al$).

\paragraph{Step 3.}
{\em There is a constant $\De_3>0$ and an integer $M_2=M_2(\al)$
such that if $J\subset I$ is an interval with
$|J|\geq \De_2\sqrt{\al}$, then for  every  $\t\in S^1$}
 $$\left|\f^{M_2}(\{\t\}\times J)\right|\geq \De_3.$$

\smallskip
\noindent Arguing as in Step 2 we can prove an analog to
(\ref{e.step2})
\[
\left|\f^N(\{\t\}\times J)\right|\geq \frac{\De_2}{4}
\sqrt\al\;\al^{-1+\eta}|J|\geq \frac{\De_2}{4} \al^{\eta}.
\]
Letting $R\geq 1$ be the first integer for which
$\f^{N+R}(\{\t\}\times J) $ intersects $[-\sqrt\al,\sqrt\al]$
we have
\[
\left|\f^{N+R} (\{\t\}\times J) \right|\geq\left\{
\begin{array}{ll}
C_2\tau^R\al^{\eta} & \mbox{if }
\f^{N+R}(\{\t\}\times J)\subset\left(S^1\times[-\de,\de]\right),\\
\de-\al^{\eta} &\mbox{otherwise.}
\end{array}
\right.
\]
In both cases we have that the $x$-component of
$\f^{N+R}(\{\t\}\times J)$ contains  some interval $L$ not
intersecting $(-\sqrt\al,\sqrt\al)$ with an end point at
$-\sqrt\al$ or $\sqrt\al$ and whose length is at least $C_2\al^\eta$.

From now on we  use $C$ to denote any large constant
depending only on the map $q$. Take $l\geq 1$ the smallest integer for
which $z=q^l(0)$ is a periodic point for $q$ and let $k\geq 1$
be its period. Denote $\rho^k=|(q^k)^\prime(z)|$ and note
that by \cite{S} we must have $\rho>1$. Fix $\rho_1,\rho_2>0$ with
$\rho_1<\rho<\rho_2$ and $\rho_1>\rho_2^{1-\eta/2}$,
and take $\de_0>0$ small enough in order to obtain
$$\rho_1^{k}<\prod_{j=0}^{k-1}\left|\partial_xf(\f^i(\sigma,y))
\right|< \rho_2^{k},\quad\mbox{whenever}\quad|y-z|<\de_0$$
(and $\al $ sufficiently small).
Since 0 is pre-periodic for $q$, there exists
some constant $\e >0$ such that $|q^j(0)|>\e$ for every $j>0$.
From this we deduce
\begin{equation}
|x_1|,\dots,|x_{l-1}|>\frac{\e}{2},\quad\mbox{whenever}\quad
x\in L,\label{xi}
\end{equation}
as long as $\al $ is sufficiently small.
By (\ref{assumption}) and (\ref{alfa})
we may write $\partial_x f(\t,x)=x\psi(\t,x)$
with $|\psi+2|<\al $ at every point $(\t,x)\in S^1\times I$.
This, together with (\ref{xi}), gives for every $x\in L$
\begin{equation}
\prod_{j=0}^{l-1}|\partial_xf(\t_j,x_j)|\geq \frac{1}{C}|x|,\label{app}
\end{equation}
and so we have for some $x\in L$
$$
\left|\f^{N+R+l}(\{\t\}\times J)\right|
\geq \prod_{j=0}^{l-1}|\partial_xf(\t_{N+R+j},x_j)|\cdot |L|
\geq \frac{1}{C}|x| \al^{\eta}\geq \frac{1}{C}\al^{1/2+\eta}.
$$
For $(\t,x)\in S^1\times I$ and $i\geq 0$ we
denote $d_i=|x_{l+ki}-z|.$ Take $\de_1>0$ and $\al $ sufficiently
small in such a way that
$$|x|<\de_1\;\Rightarrow\; d_0\leq Cx^2+C\al <\de_0.$$
If $(\t,x)$ and $i\geq 1$ are such that $|x|<\de_1$ and
$d_0,\dots,d_{i-1}<\de_0$, then
$d_i\leq \rho_2^kd_{i-1}+C\al $
and so, inductively,
$$
d_i\leq (1+\rho_2^k+\cdots+\rho_2^{k(i-1)})C\al+\rho_2^{ki}d_0
\leq \rho_2^{ki}(C\al +Cx^2).
$$
In  particular for the points $x=\pm\sqrt\al$ we have
$d_i\leq \rho_2^{ki}C\al $. Now we take $N_0={N_0}(\al)\geq 1$
the smallest integer for which $\rho_2^{k{N_0}(\al)}C\al \geq \de_0/2$.
This choice of ${N_0}$ implies
\begin{equation}
d_i<\de_0/2 \quad\mbox{for}\quad i=0,\dots,N_0
-1.\label{di}
\end{equation}
Now we consider the following two possible
cases:
\begin{enumerate}
{\em \item  $\f^{l+ki}(\{\t_{N+R}\}\times L)\subset
(z-\de_0,z+\de_0)$ for every $i\in\{0,\cdots,{N}_0-1\}.$}

This implies that
\begin{eqnarray*}
\left| \f^{l+k N_0}(\{\t_{N+R}\}\times L)\right|
 &\geq&\rho_1^{k {N}_0}\left|\f^{l}(\{\t_{N+R+k {N}_0}\}\times
    L)\right|\\
 &\geq& \rho_2^{(1-\eta/2)kN_0}\frac{1}{C}\al^{1/2+\eta}\\ &\geq&
\frac{1}{C}\alpha^{-1+\eta/2}\frac{1}{C}\al^{1/2+\eta}\\
&\geq& \frac{1}{C}\alpha^{-1/2+3\eta/2}\\ &\gg &\frac{1}{C}
\end{eqnarray*}
(recall that $\eta<1/3$).

{\em \item  $\f^{l+ki}(\{\t_{N+R}\}\times L)\not\subset
(z-\de_0,z+\de_0)$ for some $i\in\{0,\cdots,{N}_0-1\}.$}

Since $d_i\leq \de_0/2$,
it follows that
$$
\left| \f^{l+ki}(\{\t_{N+R}\}\times L)\right|\geq\de_0-\de_0/2=\de_0/2.
$$
\end{enumerate}
In both cases we have some integer $N_1\leq N_0$ for which
$$\left| \f^{l+kN_1}(\{\t_{N+R}\}\times L)\right|\geq
\frac{1}{C}.$$ Thus, taking $M_2=N+R+l+kN_1$ we have
$$\left| \f^{M_2}(\{\t\}\times J)\right|\geq \frac{1}{C}.$$

\paragraph{Step 4.}
{\em There is an integer $M_3=M_3(\al)$ such that if $J\subset I$
is an interval with $|J|\geq \De_3$, then for every  $\t\in S^1$}
$$
\left|\f^{M_3}(\{\t\}\times J)\right|=
(\{\t_{M_3}\}\times I)\cap \La.
$$

\smallskip
\noindent Since we are taking $a_0$
a Misiurewicz parameter, it follows that the pre-orbit of the
repelling fixed point $P$ of $q$ is dense. So,
there is some integer $n_1(\De_3)\geq 1$ such that for
every interval $J\subset I$ with $|J|\geq\De_3$
we have that $q^{n_1(\De_3)}(J)$ covers a neighbourhood of $P$ with a
definite size (depending only on $\Delta_3$). By a finite number
of iterates $n_2(\Delta_3)$ we transform this neighbourhood
in the whole interval $q^2(I)$. Hence, taking
$M_3=n_1(\De_3)+n_2(\Delta_3)+1$ we have by continuity
$\f^{M_3}(\{\t\}\times J)=(\{\t_{M_3}\}\times I)\cap\Lambda$
for sufficiently small $\al$.

\medskip
Now it suffices to take
 $M(\al)=M_1(\al)+M_2(\al)+M_3(\al)$ and we
complete the proof of Proposition \ref{p.whole}.
\cqd

Now we are in conditions to prove that the maps $\f\in\cn$
are topologically mixing.
Let $A$ be an open set in $S^1\times I$.
Since the partitions $\cs_n$ have diameters
converging to zero when $n$ goes to infinity, there
must be some $n\geq 1$ and $S\in\cs_n$ for which $S\subset
A$. Hence, taking $n(A)=h_n(S)+M$ ($M$ given by Proposition
\ref{p.whole}) it follows from Proposition \ref{p.whole}
that $\f^{n(A)}(A)=\Lambda$.


\section{Ergodicity}

In this section we prove the ergodicity of the maps
$\f\in\cn $ with respect to Lebesgue measure. We start by
proving some auxiliary results.

\cle\label{suf}
Let $B$ be a Borel subset of $S^1\times I$ such that
$\f^{-1}(B)=B$.
\begin{enumerate}
\item If $m(B\cap\Lambda)=0$ then $m(B)=0$.
\item If $\f^n(R)=\Lambda$ for some $n\geq 1$ and
$R\subset S^1\times I$, then $B\cap\Lambda\subset\f^n(B\cap R)$.
\end{enumerate}
\fle
 \dem  Since we have $\Lambda=\f^2(S^1\times I)$, it
follows that
$$A=\f^{-2}(A)=\f^{-2}\left(A\cap\f^2(S^1\times I)\right)
=\f^{-2}\left(A\cap\Lambda\right).$$ Thus, if
$m(A\cap\Lambda)=0$, then $m(A)=0$, and so we have proved
the first item.

Now let $x\in B\cap\Lambda$. Since $\f^n(R)=\Lambda$, there
must be some $z\in R$ for which $\f^n(z)=x$. On the other
hand, since $\f^{-1}(B)=B$ we have $\f^{-n}(B)=B$, and so
$z$ belongs to $B$. Hence $x\in\f^{n}(B\cap B)$.
\cqd

Now we prove a general result that will play an important
role in the proof of the ergodicity with respect to
Lebesgue measure.

\cpr \label{aprox}
Let $X$ be a metric space, $\mu$ be a Borel measure on $X$,
and $\cp=\{P_1,\dots,P_r\}$ be a partition of $X$ into Borel subsets.
Assume that $(\cs_n)_{n\geq 1}$ are partitions of $X$
such that diam$\:(\cs_n)\rightarrow 0$ when $n\rightarrow\infty$.
Then, for each $n\geq 1$ there is a partition
$\{Q^n_1,\dots,Q^n_r\}$ of $X$ such that for $i=1,\dots,r$
\begin{enumerate}
\item $Q^n_i$ is a union of atoms of $\cs_n$.
\item $\lim_{n\rightarrow\infty}\mu(Q^n_i\triangle P_i)=0$.
\end{enumerate}
\fpr

\dem
Take an arbitrary $\epsilon>0$. Since $\mu$ is a regular
measure, there are compact sets $K_1,\dots,K_r\subset X$ with
$$K_i\subset P_i\quad\mbox{and}\quad\mu(P_i\setminus K_i)<\epsilon$$
 for $i=1,\dots,m.$
Let
$$\delta=\inf_{i\neq j}d(K_i,K_j)>0$$
and take $n_0\geq 1$ such that
$$\mbox{diam}(\cs_n)<\delta/2\quad \mbox{for}\quad n\geq n_0.$$
 For $n\geq n_0$ we divide $\cs_n$ into $r$ groups,
whose unions we call $Q^n_1,\dots,Q^n_r$, by putting
$S\subset Q^n_i$ if $S\in\cs_n$ intersects $K_i$.
Note that each $S\in\cs_n$ intersects at most one $K_i$. If
it does not intersect any $K_i$, then we include
it arbitrarily in some $Q_i^n$. We have
\begin{eqnarray*}
\mu(Q^n_i\triangle P_i)
&=&\mu(Q^n_i\setminus P_i)+\mu(P_i\setminus Q^n_i)\\
&\leq& \mu(X\setminus \cup_{i=1}^{r}K_i)+\mu(P_i\setminus K_i)\\
&\leq& (r+1)\epsilon.
\end{eqnarray*}
Since $\epsilon>0$ is arbitrary and $r$ is fixed, we have
proved the result.
\cqd

The following corollary may be thought as a similar result to
the Lebesgue density theorem for balls in Euclidean spaces.

\cco\label{c.dens}
Let $B$ be a Borel subset of $X$ with $\mu(B)>0$.
Then, for every $\epsilon>0$ there is an integer $n_\epsilon\geq 1$
such that for each $n\geq n_\epsilon$ there is $S\in\cs_n$ with
$$\mu(B^c\cap S)<\epsilon\mu(S).$$
\fco
\dem
Assume by contradiction that there are $\e_0>0$
and a sequence of integers $(n_k)_{k\geq 1}$ going to $+\infty$ such
that
\begin{equation}
\forall\,k\geq 1\;\forall\,S_{n_k}\in\cs_{n_k}
\;:\;\mu(B^c\cap S_{n_k})\geq\e_0\mu(S_{n_k}).\label{a1}
\end{equation}
We know from Lemma \ref{aprox} that for every $k\geq 1$ there
is a partition $\{Q^{n_k}_1,Q^{n_k}_2\}$ of $X$ such that
$Q^{n_k}_1$, $Q^{n_k}_2$ are unions of atoms of $\cs_{n_k}$,
and
$$
\lim_{k\rightarrow +\infty}\mu(Q^{n_k}_1\triangle B)=0,\quad
\lim_{k\rightarrow +\infty}\mu(Q^{n_k}_2\triangle B^c)=0.
$$
Take $\epsilon=\e_0\mu(B)/(1+\epsilon_0)$ and $k_0\geq 1$
sufficiently large in order to
\begin{equation}
\mu(Q^{n_{k_0}}_1\triangle B)<\e
\quad\mbox{and}\quad
\mu(Q^{n_{k_0}}_2\triangle B^c)=\e.
\label{a3}
\end{equation}
Letting $\cs_{n_{k_0}}=\{S_i\}_{i\in\NN}$, we know that there is some
${\mathbb I}\subset \NN$ for which
$$
Q^{n_{k_0}}_1=\bigcup_{i\in {\mathbb I}}S_i\quad\mbox{and}\quad
Q^{n_{k_0}}_2=\bigcup_{i\in \NN\,\setminus {\mathbb I}}S_i.
$$
From (\ref{a1}) we have in particular
$$
\mu(B^c\cap S_{i})\geq\e_0\mu(S_{i})
$$
for every $i\in {\mathbb I}$, and so summing over all $i\in {\mathbb I}$
we have
\begin{equation}
\mu(B^c\cap Q_1^{n_{k_0}})\geq\e_0\mu(Q_1^{n_{k_0}}).\label{a4}
\end{equation}
Finally, from (\ref{a3}) and (\ref{a4}) we get
$$\e>\mu(B^c\cap Q_1^{n_{k_0}})\geq\e_0\mu(Q_1^{n_{k_0}})
\geq \e_0\left(\mu(B)-\mu(B\setminus Q_1^{n_{k_0}})\right)>\e_0(\mu(B)-\e),$$
which gives a contradiction (recall our choice of $\e$).
\cqd

It will be useful to have the following distortion result,
whose role is essentialy to state the non dependence
on the partition $\cs_n$
of the constant in Proposition \ref{p.dist0} .

\cpr \label{p.dist}
There is some constant $\De>1$ such
that for every $n\geq 1$, $S\in \cs_n$ and
$(\t,x),(\si,y)\in S$ we have
$$
\frac{1}{\De}\leq
\left|\frac{J_n(\t,x)}{J_n(\si,y)}\right|\leq \De,
$$
where $J_n$ is the Jacobian of $\f^{h_n(S)}|S$.
\fpr
\dem
We observe that the constant in \cite[Proposition
4.2]{A} that bounds the distortion does not depend on the
integer $p$ that we have used for starting the construction of
the partition (it essentially depends on the expansion
rates of the maps $g$ and $f$). This means that the same
proof of Proposition \ref{p.dist0} applies to this
situation with the constant $\De>0$ not depending on $n$.
\cqd

Now we are in conditions to prove the ergodicity of
the maps $\f\in\cn$ with respect to the Lebesgue measure.
Let $B$ be a Borel subset of $S^1\times I$ with
$\f^{-1}(B)=B$ and having positive Lebesgue measure. We
need to prove that the Lebesgue measure of $B^c=(S^1\times
I)\setminus B$ is equal to zero. From the first item of
Lemma \ref{suf} it
suffices to prove that $m(B^c\cap \Lambda)=0$. Take any
$\e>0$ small. It follows from Corollary \ref{c.dens} that
there are $n_\e\geq 1$ and $S\in\cs_{n_\e}$ for which
$$
m(B^c\cap S)<\e m(S).
$$
Letting $h=h_{n_\e}(S)$ we have
from Proposition \ref{p.whole} that $\f^{h+M}(S)=\Lambda$.
Thus, applying the second item of Lemma \ref{suf} we obtain
$$
m(B^c\cap\Lambda)\leq
m\left(\f^{h+M}(B^c\cap S)\right).
$$
Fixing some $(\t_0,x_0)\in S$, we deduce from Proposition
\ref{p.dist}
\begin{eqnarray*}
m\left(\f^{h}(B^c\cap S)\right)
&=& \int_{B^c\cap S}\left|J_{n_\e}(\t,x)\right|dm(\t,x)\\
&\leq &\De\int_{B^c\cap S}\left|J_{n_\e}(\t_0,x_0)\right|dm(\t,x)\\
&\leq &\De\left|J_{n_\e}(\t_0,x_0)\right|m\left(B^c\cap S\right).
\end{eqnarray*}
Similarly we prove that $$ m\left(\f^{h}( S)\right) \geq
\frac{1}{\De}\left|J_{n_\e}(\t_0,x_0)\right|m\left(
S\right). $$ Hence $$m\left(\f^h(B^c\cap S)\right)\leq
\frac{m\left(\f^h(B^c\cap S)\right)}{m\left(\f^h(
S)\right)}\leq \De^2 \frac{m\left(B^c\cap S\right)}{m\left(
S\right)}\leq \De^2\e ,$$ which finally gives
$$m\left(\f^{h+M}(B^c\cap S)\right)\leq (d+\al)^M4^M
m\left(\f^h(B^c\cap S)\right)\leq (d+\al)^M4^M \De^2\e.$$
Since this holds for $\e $ arbitrarily small and $M$ is
fixed, the proof is complete.

\bigskip

\noindent Jos\'e Ferreira Alves\\
Centro de Matem\'atica da Universidade do Porto\\
Pra\c ca Gomes Teixeira,
4099-002 Porto, Portugal\\
{\tt jfalves@fc.up.pt}

\bigskip

\noindent Marcelo Viana\\ IMPA, Estrada Dona Castorina,
110\\ 22460-320 Rio de Janeiro, Brasil\\ {\tt
viana@impa.br}


\begin{thebibliography}{ABV}

\bibitem[A]{A} J. F. Alves, {\em SRB measures for non-hyperbolic
systems with multidimensional expansion}, to appear in Ann.
Sci. de l'ENS.

\bibitem[ABV]{ABV} J. F. Alves, C. Bonatti, M. Viana, {\em SRB measures
for partially hyperbolic systems with mostly expanding central direction},
to appear in Invent. Math.

\bibitem[AP]{AP} A. Andronov and L. Pontryagin, {\em Syst\`emes grossiers},
Dokl. Akad. Nauk. USSR, 14 (1937), 247--251.

\bibitem[BV]{BV} C. Bonatti, M. Viana, {\em SRB measures for partially
hyperbolic systems with mostly contracting central direction}, to appear
Israel J. Math.

\bibitem[D]{D} D. Dolgopyat, {\em On dynamics of mostly contracting diffeomorphisms}, preprint 1998.

\bibitem[G]{G} E. Giusti, {\em Minimal surfaces and functions of bounded
variation}, Birk\"auser Verlag, Basel-Boston, Mass., 1984.

\bibitem[HPS]{HPS} M. Hirsch, C. Pugh, M. Shub, {\em Invariant manifolds},
Lect. Notes in Math. 583, Springer Verlag, 1977.

\bibitem[LY]{LY} A. Lasota and J.A. Yorke, {\em On the existence of invariant
measures for piecewise monotonic maps}, Trans. Amer. Math. Soc. {\bf 186}
(1973), 481-488.

\bibitem[PS]{PS} J. Palis and S. Smale,
{\em Structural stability theorems}, in {\em Global Analysis},
Proc. Sympos. Pure Math. XIV (Berkeley 1968), Amer. Math. Soc., 223--232, 1970.

\bibitem[PT]{PT} J. Palis and F. Takens,
{\em Hyperbolicity and sensitive-chaotic dynamics at homoclinic bifurcations},
Cambridge University Press, 1993.


\bibitem[S]{S} D. Singer, {\em Stable orbits and bifurcations of maps
of the interval}, SIAM J. Appl. Math. {\bf 35} (1978), 260-267.

\bibitem[V1]{V} M. Viana, {\em Multidimensional nonhyperbolic
attractors}, Publ. Math. IHES {\bf 85} (1997), 63-96.

\bibitem[V2]{V2} M. Viana, {\em Stochastic dynamics of deterministic
systems}, Lect. Notes XXI Braz. Math Colloq., IMPA, 1997.

\end{thebibliography}
\end{document}